\documentclass[12pt,a4paper]{article}
\usepackage{amsmath, amssymb,  latexsym,amsthm,amsfonts}
\usepackage{graphicx}

\newtheorem{theorem}{Theorem}[section]
\newtheorem{lemma}[theorem]{Lemma}
\newtheorem{proposition}[theorem]{Proposition}

\newtheorem{assumption}[theorem]{Assumption}
\newtheorem{conjecture}[theorem]{Conjecture}
\newtheorem{corollary}[theorem]{Corollary}

\newtheorem{question}[theorem]{Question}
\newtheorem{remark}[theorem]{Remark}
\newtheorem{remarks}[theorem]{Remarks}

\allowdisplaybreaks

\def\neweq#1{\begin{equation}\label{#1}}
\def\endeq{\end{equation}}

\def \Epsilon{{\mathcal E}}

\def\omp {{\omega_\perp}}
\def\dep{{\delta_{\perp}}}

\def\ep {{\epsilon}}

\def \rz { {\mathbb R}}

\def\Og {{\cal O}} 

\def \rz {{\mathbb R}}

\def \div{{\rm div\,}}
\def \cb {{\bf c}}

\def\neweq{\begin{equation}}
\def\endeq{\end{equation}}
\def\beq{\begin{equation}}
\def\eeq{\end{equation}}

\newcommand {\ar}{\rightarrow}
\newcommand {\pa}{\partial}

\numberwithin{equation}{section}
\begin{document}
 \title{\bf On mathematical models for
Bose-Einstein
  condensates
 in optical lattices (expanded version)}
 \author{Amandine Aftalion\footnote{CNRS,  CMAP, Ecole Polytechnique, F-91128 Palaiseau
 cedex, France} \ and Bernard Helffer\footnote{Laboratoire de
     Math\'ematiques,
 Univ Paris-Sud et CNRS, Bat 425. 91 405 Orsay Cedex France.}}

 \maketitle

\begin{abstract}
Our aim is to analyze the various energy functionals appearing in
the physics literature and describing the behavior of a
Bose-Einstein condensate
 in an optical lattice. We want to justify the use of some reduced models.
For that purpose, we will use the semi-classical analysis developed
for linear problems related to the Schr\"odinger operator
 with periodic potential or multiple wells potentials.
  We  justify, in some asymptotic regimes,
   the reduction to low dimensional problems and analyze the reduced problems.

\end{abstract}

\tableofcontents

\section{Introduction}\subsection{The physical motivation for Bose-Einstein condensates in optical lattices}
Superfluidity and superconductivity are two spectacular
manifestations of quantum mechanics at the macroscopic scale. Among
their striking characteristics is the existence of vortices with
quantized circulation. The physics of such vortices is of tremendous
importance in the field of quantum fluids and extends beyond
 condensed matter physics. The advantage of ultracold gaseous Bose-Einstein condensates
is to allow tests in the laboratory to study various aspects of
macroscopic quantum physics.
 There is a large body of research, both experimental,
  theoretical and mathematical on vortices in Bose-Einstein
condensates \cite{PeSm,PiSt,Af,LSSY}. Current physical interest is
in the investigation of very small atomic assemblies, for which one
would have one vortex per particle,  which is a challenge in terms
of detection and signal analysis. An appealing option consists in
parallelizing the study, by producing simultaneously a large number
of micro-BECs rotating at the various nodes of an optical lattice
\cite{Sn}. Experiments are under way. A major topic is the
transition from a Mott insulator phase to a superfluid phase. We
refer to the paper of Zwerger \cite{Z} and the references therein
 for more details. Our framework of study will be in the mean field
  regime where the condensate can be described by a Gross Pitaevskii
  type energy with a term modeling the optical lattice potential.
   The  mean field description of
a condensate by the Gross Pitaevskii energy has been derived as the
limit of the hamiltonian for $N$ bosons, when $N$ tends to infinity
\cite{LSY,LS} in the case of a condensate without optical
 lattice. The scattering
 length $a_N$ of the interaction in the $N$-body problem is such that
 $Na_N \rightarrow g$. The rigorous derivation in the case of an
 optical lattice where there are fewer atoms per site is
 nevertheless open.
 In a
one-dimensional optical lattice, the condensate splits into a stack
of weakly-coupled disk-shaped condensates, which leads to some
intriguing analogues with high-Tc superconductors due to their
similar layered structure \cite{SnSt,SnSt1,KMPS,ABB1,ABB2,ABS}. Our
aim, in this paper, is to address
 mathematical models that describe a
 BEC in an optical lattice.
 The theory which we will develop is inspired by  a series of physics papers \cite{Sn,SnSt,SnSt1,KMPS,STKB}. We want
 to justify their reduction to simpler
  energy functionals in certain regimes of parameters and in
  particular understand the ground state energy.

The ground state energy of a rotating Bose-Einstein condensate is given by
the minimization of
\begin{multline}\label{energyground}
Q_{\Omega}(\Psi):= \\
\int_{\mathbb R^3}
\left(\frac{1}{2} |\nabla \Psi-i\Omega\times{\bf r}\Psi |^2-\frac12 \Omega^2
r^2\,|\Psi|^2\, + (V({\bf r})+W_\epsilon(z)) |\Psi|^2 + g |\Psi|^4\right)\,dxdydz\,,
 \end{multline}
under the constraint
\beq \int_{\mathbb R^3} |\Psi(x,y,z)|^2\,dxdydz =1\,,
\eeq
where
\begin{itemize}
\item
$
r^2 = x^2+y^2\,, \,{\bf r} = (x,y,z)\,,
$
\item
 $\Omega \geq 0$ is the rotational
velocity along the $z$ axis,
\item
$\Omega\times {\bf r}= \Omega (-y,x,0)\,,$
\item  $g\geq 0 $ is the scattering length.
\end{itemize}The experimental device leading to the realization of
optical lattices requires a trapping potential 
  $V({\bf r})$
given by
\begin{equation}\label{1.1}
V({\bf r}) = \frac 12 \left(\omega_{\perp}^2 r^2 + \omega_z^2
z^2\right),
\end{equation}
corresponding to the magnetic trap. We assume that
 the radial trapping frequency  is much larger than the axial
trapping frequency, that is
\begin{equation}\label{1.2}
0\leq \omega_z <<  \omega_{\perp}\,.
\end{equation} We will always assume the condition~:
\begin{equation}\label{conditionsurOmega}
0 \leq \Omega < \omp
\end{equation} for the existence of a minimizer. The trapping has to
be stronger than the centrifugal force. The presence of the one
dimensional optical lattice in the $z$ direction is modeled by
\begin{equation}\label{1.4}
 W_\epsilon (z) = \frac{1}{\epsilon^2} {\bf w} (z)\,,
\end{equation}
where $\frac{1}{\epsilon^2}$ is the lattice depth\footnote{called
$V_z$
 in Snoek \cite{Sn}}, and $w$ is  a positive
 $T$-periodic function.
 In  the
whole paper, we will assume~:
\begin{assumption}\label{assw}~\\
The potential ${\bf w}$ is a $C^\infty$ even, non negative  function on $\mathbb R$ which is
 $T$-periodic and admits as unique minima the points $k T$ ($k\in
 \mathbb Z$). Moreover these minima are non degenerate. Thus,
\begin{equation}\label{1.5}
{\bf w}(z+T)={\bf w}(z)\,,\; {\bf w}(0)=0\,,\; {\bf w}''(0) >0\,,\; {\bf w}(z) >0 \mbox{ if }
z\not\in T\mathbb Z\,.
\end{equation}
\end{assumption}
An example is
\begin{equation}
{\bf w}(z) =\sin^2(\frac{2\pi z}{\lambda})
\end{equation}
and $\lambda$ is the wavelength of the laser light.
 The optical potential $W_\epsilon$
creates a one-dimensional lattice of wells separated by
 a distance  $ T = \lambda/2\,.$
We will assume that $\epsilon$ tends to $0$
 (this means deep lattice) and that $T$ is fixed.
  Furthermore, we assume that
the lattice is deep enough so that it dominates over the magnetic
trapping potential in the $z$ direction and that the number of sites
is large. Thus we will, in this paper, ignore the magnetic trap in
the $z$ direction, and this
 will correspond to the case
\begin{equation}\label{omegaznul}
\omega_z =0\,.
\end{equation}
We will
 mainly discuss, instead of a problem in $\mathbb R^3$,  a periodic problem
 in the $z$ direction, that is in $\mathbb R_{x,y}^2 \times [-\frac
 T2,\frac T2[$,
 where $T$ corresponds
 to the period of the optical lattice, or in $\mathbb R_{x,y}^2 \times
 [-\frac{NT}{2},\frac{NT}{2}[$ for  a fixed integer $N\geq 1$.
Therefore, we  focus (see however Subsection \ref{ss3D} for a
 justification of this choice) on the minimization of the functional
\begin{multline}\label{defQBEperN}
Q_{\Omega}^{per,N}(\Psi):= \\ \int_{\mathbb R^2_{x,y}\times
  ]-\frac{NT}{2}, \frac{NT}{2}[}
\left( \frac{1}{2} |\nabla \Psi-i\Omega\times {\bf r}\Psi|^2-\frac12 \Omega^2
r^2 |\Psi|^2 + (V({\bf r})+W_\epsilon(z)) |\Psi|^2 + g |\Psi|^4\right)\,dxdydz\,,
\end{multline}
under the constraint
 \beq
 \int_{\mathbb R^2_{x,y}\times  ]-\frac{NT}{2}, \frac{NT}{2}[}
|\Psi(x,y,z)|^2\,dxdydz =1\,,
\eeq with
\beq\label{Vp}V({\bf r})=\frac 12
\omp^2r^2\,,\eeq
the potential  $W_\ep$  given by
 (\ref{1.4})-(\ref{1.5}), and the wave function $\Psi$ satisfying
\begin{equation}\label{propper}
\Psi (x,y, z+NT)= \Psi (x,y,z)\,.
\end{equation}
This functional has  a minimizer in the unit sphere of  its natural form domain (see
 \eqref{Operunit}
 for its description)
$\mathcal S_{\Omega}^{per,N}$ and we call
\begin{equation}\label{defEperN}
E^{per,N}_\Omega = \inf_{\Psi \in \mathcal S_{\Omega}^{per, N}}
Q_{\Omega}^{per,N}(\Psi)\,.
\end{equation}
\paragraph{Notation}~\\
In the case $N=1$, we will write more simply
\beq\label{defQBEper}
Q_{\Omega}^{per}:=Q_{\Omega}^{per,(N=1)}\,,\; E_\Omega^{per} :=E_\Omega^{per,(N=1)}\,.
\eeq
When $\Omega=0$, we will sometimes omit the reference to $\Omega$.\\

Our aim is to justify that the ground state energy can be well approximated
 by the study of simpler models introduced
in physics papers \cite{Sn,SnSt,KMPS} and to measure the error which
 is done in the approxiamtion.\\
 For that purpose, we
 will describe how, in certain regimes, the semi-classical analysis
developed for linear problems related to the Schr\"odinger operator
 with periodic potential or multiple wells potentials is relevant:
  Outassourt \cite{Ou}, Helffer-Sj\"ostrand \cite{He,DiSj}
 or for an alternative approach \cite{Si}.

\subsection{The linear model}

The linear model which appears naturally is a selfadjoint realization
associated with the differential operator~:
\begin{equation}\label{3.1}
H_\Omega = H^\Omega_{\perp} + H_z\,,
\end{equation}
with
\begin{equation}\label{3.2}
H^\Omega_{\perp}:= - \frac{1}{2} \Delta_{x,y} + \frac 12 \omega_\perp^2 r^2
- \Omega L_z \,,
\end{equation}
\begin{equation}
L_z =  i   (x\partial_y - y \partial_x)\,,
\end{equation}
and
\begin{equation}\label{3.4}
H_z:=  - \frac{1}{2} \frac{d^2}{dz^2} + W_\epsilon(z)\,.
\end{equation}
In the transverse direction, we will consider the unique natural
selfadjoint extension in
 $L^2(\mathbb R_{x,y}^2)$
 of the positive operator  $H^\Omega_{\perp}$ by keeping the same
 notation. In the longitudinal direction,
 we will consider specific realizations of $H_z$ and in particular the
 $T$-periodic problem or more generally the $(NT)$-periodic
 problem attached to $H_z$ which will be denoted by $H_z^{per}$ and
 $H_{z}^{per,N}$ and we keep the notation $H_z$ for the problem on the
 whole line.\\
So our model will be the self-adjoint operator
\begin{equation}
H^{per,N}_\Omega = H_\perp^\Omega + H_z^{per,N}\,.
\end{equation}

In this situation with separate variables, we can split the spectral
analysis, the spectrum
of $H^{per,N}_\Omega$ being the closed set
\begin{equation}\label{3.6}
\sigma (H^{per,N}_\Omega):= \sigma (H^\Omega_{\perp}) + \sigma(H_z^{per,N})\,.
\end{equation}

The first operator $H^\Omega_{\perp}$ is a harmonic oscillator with
discrete spectrum as we will explain in Section \ref{s2}. Under
Condition \eqref{conditionsurOmega}, the bottom of
its spectrum
 is given by
\begin{equation}
\lambda_1^{\perp} :=\inf (\sigma(H^\Omega_{\perp} )) =
\omega_\perp \,,
\end{equation}
hence is independent of $\Omega$.\\
A corresponding groundstate is the Gaussian $\psi_{\perp} =
\left(\frac{\omega_\perp}{\pi}\right)^\frac 12 \exp - \frac
{\omega_\perp}{2} r^2$. The gap between the ground state energy and
the second eigenvalue (which has multiplicity $1$ or $2$) is given by
\begin{equation}
\delta_\perp:=\lambda_{2,\Omega}^\perp - \lambda_1^\perp =
\omega_\perp-\Omega \,.
\end{equation}

The properties of  the periodic Hamiltonian $H_z^{per,N}$, which will be
described in Subsection~\ref{ssha} (Formulas \eqref{ha5} and
\eqref{longgap}
 for the physical model), depend on the value of $N$.
In the
case $N=1$, we call the groundstate of $H_z^{per}$ $\phi_1(z)$ and the ground
energy (or lowest eigenvalue)
 $\lambda_{1,z}$.
  In the semi-classical regime $\epsilon \to 0$, $\lambda_{1,z}$ satisfies
\beq \label{asymp1} \lambda_{1,z}\sim \frac{c}{\epsilon}, \eeq for
some $c>0$. The splitting  $\delta_z$ between the groundstate energy
 and the first excited eigenvalue satisfies
\beq\label{asymp2} \delta_z \sim \frac {\tilde c}{\epsilon}\,, \eeq
for some $\tilde c >0$.\\
  For $N>1$, the groundstate energy of $H_z^{per,N}$ is unchanged and the
  corresponding groundstate $\phi_1^N$ is the periodic extension of
  $\phi_1$ considered as an $(NT)$-periodic function. More precisely,
 in order to have the $L^2$- normalizations, the relation is
\begin{equation}\label{deunaN}
\phi_1^N = \frac{1}{\sqrt{N}} \phi_1\,,
\end{equation}
on the line. But we have now $N$ exponentially close eigenvalues of
the order of $\lambda_{1,z}$
 lying in the first band of the spectrum of
 $H_z$  on
 the whole line. They are separated from the $(N+1)$-th by a splitting
 $\delta_z^{N}$ which satisfies~:
\beq \label{defphiN} \delta_z^{N} = \delta_z + \widetilde \Og (\exp
- S/\ep))\,. \eeq Here the notation $\widetilde \Og (\exp - S/\ep))$
means \beq\label{defOtilde} \widetilde \Og (\exp - S/\ep)) =
\Og(\exp -S'/\ep) \,,\; \forall S'<S\,. \eeq
 The
 first $N$ eigenfunctions
 satisfy
\beq
\phi_\ell^N(z+T)=\exp (\frac{2i\pi (\ell-1)}{N})\;\phi_\ell^N(z)\,,\;
\mbox{ for }
\ell =1,\dots, N \,,\eeq
corresponding to the special values $k= \frac{2\pi(\ell -1)}{NT}$ of what will be called
 later a $k$-Floquet condition (see \eqref{A.2}).

We will sometimes use another real orthonormal basis (called
$(NT)$-periodic Wannier functions basis) $(\psi_j^N)$
($j=0,\dots,N-1$)  of the spectral space attached to the first $N$
eigenvalues. Each of these $(NT)$-periodic functions have the
advantage to be localized (as $\epsilon \ar 0$) in a specific well
of $W_\epsilon$ considered as defined on $\mathbb R/(NT)\mathbb Z$.

\subsection{The reduced functionals}
We want to prove the reduction to lower dimensional functionals
  by using the spectral analysis of the linear problem. There are two natural ideas to
   compute upper bounds, and thus find these functionals. We can
\begin{itemize}
\item  either
    use test functions  of the type
\beq\label{tfa}\Psi(x,y,z) = \phi(z)
 \psi_{\perp}(x,y)\,,\eeq where $\psi_\perp$ is the first normalized
 eigenfunction
  of $H_\perp^\Omega$  and minimize among all possible
    $L^2$-normalized
  $\phi (z)$ to obtain a 1D-longitudinal reduced problem,
\item  or use
\begin{itemize}\item in the case $N=1$,
 \beq\label{tfb}\Psi(x,y,z) = \phi_1(z)
 \psi (x,y)\eeq  where $\phi_1$ is the first eigenfunction
  of $H_z^{per}$ and minimize among all possible $L^2$-normalized
  $\psi (x,y)$ to obtain a 2D transverse reduced problem,
\item or in the case
  $N\geq 1$ \beq\label{tfbN}\Psi(x,y,z) = \sum_{j=0}^{N-1} \psi_j^N(z)
 \psi_{j,\perp} (x,y)\eeq where $\psi_j^N(z)$
  is the orthonormal basis of Wannier functions mentioned above,
  and minimize on the suitably normalized
 $\psi_{j,\perp}$'s which provide
 $N$ coupled problems. We denote by $\Pi_N$ the projection on this
  space. For $\Psi \in L^2(\mathbb R^2\times ]-\frac{NT}{2},\frac{NT}{2}[)\,$,
  we have
\beq\label{defPiN}
\Pi_N \Psi =\sum_{j=0}^{N-1} \psi_j^N(z)
 \psi_{j,\perp} (x,y)\,,
\eeq
 with 
$$  \psi_{j,\perp} (x,y) = \int_{
  ]-\frac{NT}{2},\frac{NT}{2}[}\Psi (x,y,z)  \psi_j^N(z)\,dz\,.
$$

\end{itemize}
\end{itemize}

Computing the energy of a test function of type (\ref{tfa}), we get
\beq \label{z1} Q_{\Omega}^{per,N}(\Psi) = \omega_\perp  +
\Epsilon^N_{A}(\phi) \eeq
where $\Epsilon_A^N$ is the functional on the
$NT$-periodic functions in the
 $z$ direction, defined on $H^1(\mathbb R/NT \mathbb Z)$ by
\beq\label{EpA} \phi\mapsto \Epsilon_{A}^N (\phi)= \int_{-\frac{NT}{2}}^{\frac{NT}{2}} \left(
\frac{1}{2} |\phi'(z)|^2 +
 W_\epsilon(z)|\phi(z)|^2 +  \widehat g \, |\phi(z)|^4 \right) \; dz \;
\end{equation} with
\beq \label{defwidehatg} \widehat g:= g  \left(\int_{\mathbb R^2}
|\psi_{\perp}(x,y)|^4\,
 dxdy\right) = \frac{1}{2 \pi}  g \omega_\perp.
\eeq The functional $\Epsilon_{A}^N $ is introduced by \cite{KMPS}
who analyze a particular case. Its study in the small $\ep$ limit is
one of the aims of this paper.

For test functions of type (\ref{tfb}),  we get in the case $N=1$
\beq
Q_{\Omega}^{per}(\Psi) = \lambda_{1,z}\, + \Epsilon_{B,\Omega}(\psi) \eeq  with
\begin{multline}
 \Epsilon_{B,\Omega}(\psi)\\
:= \int_{\mathbb R^2_{x,y}}\left (  \frac12|\nabla_{x,y} \psi
-i\Omega \times{\bf r}\psi|^2-\frac12 \Omega^2r^2|\psi|^2
 +\frac12 \omega_{\perp}^2 (x^2+y^2)  |\psi |^2
 + \widetilde g |\psi|^4\right )\;  dx\,
dy\,,\label{EB}
\end{multline}
and
\beq \label{defwidetildeg} \widetilde g := g \left(\int_{-\frac
    T2}^{\frac T2}
|\phi_1(z)|^4\, dz\right)\,. \eeq
 In the case $N>1$, we define $\Epsilon_{B,\Omega}^N((\psi_{j,\perp})_{j=0,\dots,N-1})$ by
\beq\label{defEBN} Q_{\Omega}^{per,N}(\Psi): = \lambda_{1,z}\;
\sum_j||\psi_{j,\perp}||^2 + \Epsilon_{B,\Omega}^N((\psi_{j,\perp}))
\eeq with \beq\label{defPsiN} \Psi = \sum_{j=0}^{N-1} \psi_j^N(z)
 \psi_{j,\perp} (x,y)\,.
\eeq Of course when minimizing over normalized $\Psi$'s, one gets
more simply the problem of minimizing \beq \label{plussimple}
Q_{\Omega}^{per,N}(\Psi) = \lambda_{1,z} +
\Epsilon_{B,\Omega}^N((\psi_{j,\perp})) \,. \eeq

As such, the energy
$\Epsilon_{B,\Omega}^N$ does not provide $N$ coupled problems but
one single energy depending on $N$ test functions. Nevertheless, in
the small $\ep$ limit, the Wannier functions are localized in each
well. Thus each function $\psi_{j,\perp}$ only interacts with its
nearest neighbors and this simplification provides $N$ coupled
problems, as suggested by Snoek \cite{Sn} on the basis of formal
computations.
 We will
 analyze their validity. This reduced functional  is somehow
 related to the Lawrence-Doniach model for superconductors (see
 \cite{ABB1,ABB2}).

\subsection{Main results}

\subsubsection{Universal  estimates and applications}
The analysis of the linear case immediately leads  to the following
trivial and universal  inequalities (which are valid for any $N$ and
any $\Omega$ such that $ 0\leq \Omega <\omega_\perp$)
\begin{equation}\label{universal}
\lambda_{1,z} + \omega_\perp
 \leq E_\Omega^{per,N} \leq  \lambda_{1,z} + \omega_\perp + I_N\,,
\end{equation}
where
\begin{equation}
I_N:=  \frac{g \omp}{2 N\pi } \left(\int_{-\frac{T}{2}}^{\frac{T}{2}}
 |\phi_1(z)|^4 dz\right)=\frac{I}{N}\,.
\end{equation}
This universal estimate is obtained by using
 the test function
$$
\Psi^{per,N} (x,y,z)=\psi_\perp(x,y) \phi_1^{N}(z)\,,
$$
where $\phi_1^{N}$ is the $N$-th normalized ground state introduced
in \eqref{deunaN} and  $\psi_\perp(x,y)$  is the  ground state of
$H_\perp^\Omega$, actually independent of $\Omega$.\\ From
(\ref{deunaN}), we have~:
\begin{equation}\label{relationphi1}
\int_{-\frac{NT}{2} }^{\frac{NT}{2}}(\phi_1^N(z))^4\, dz =
\frac{1}{N^2} \int_{-\frac{NT}{2}} ^{\frac{NT}{2}}\, \phi_1(z)^4\,
dz =\frac 1N \int_{-\frac{T}{2}} ^{\frac{T}{2}}\, \phi_1(z)^4\,
dz\,,
\end{equation}
where,  as  $\ep \ar 0$, and, under Assumption \eqref{1.5}, it can
be proved (see \eqref{normel4a}) that
\begin{equation}\label{roughphi4}
\int_{-\frac{T}{2}} ^{\frac{T}{2}}\, \phi_1(z)^4\, dz
 \sim c_4\, \ep^{-\frac 12}\,,
\end{equation}
for some explicitly computable constant $c_4>0$.\\
 Thus, we have the following estimate for $I_N$
\begin{equation}\label{asymptIN}
I_N \sim \frac{c_4}{2\pi} \, \frac{g\omp}{N}\, \ep^{-\frac 12}\,.
\end{equation}

An immediate analysis shows that $\lambda_{1,z} +
\omega_\perp$
 is a good asymptotic of $E_\Omega^{per,N}$ in the limit $\epsilon \ar
 0$ when $g$ is sufficiently small (what we can call the quasi-linear
 situation). More precisely, we have
\begin{theorem}\label{ThQL}~\\
 Under the condition that either
\begin{equation}\label{condbiza}
(QLa) \quad g << \epsilon^\frac 12\,,
\end{equation}
or
\begin{equation}\label{condbizb}
(QLb)\quad  g\omp\ep ^\frac 12 << 1\,,
\end{equation}
then we have
\begin{equation}\label{asymptQL}
E_\Omega^{per,N} = \left( \lambda_{1,z} +
\omega_\perp \right) (1+ o(1))\,,
\end{equation}
as $\epsilon$ tends to $0$.
\end{theorem}

Each of these conditions implies indeed that  $I_N$ is small relatively to
$\lambda_z$ or to $\omp$.\\

So our goal is
\begin{itemize}
\item
to have more accurate estimates than \eqref{asymptQL},
\item to analyze more  interesting cases when
 none of these two conditions is satisfied and to
 give natural sufficient conditions allowing the analysis of
 reduced models.
\end{itemize}

We are able to justify the reductions to the lower dimensional
functionals $\Epsilon_{A}^N $ and $\Epsilon_{B,\Omega}^N $ when their infimum
 is much smaller than the gap between the first two excited states
 of the linear problem in the other direction, namely in case A,
 when $m_A^N$ is much smaller than $\delta_\perp$, where
\beq\label{defmA}
m_A^N = \inf_{||\phi||=1}
\Epsilon_{A}^N (\phi)\,,
\eeq
 and in case B, when $m_{B,\Omega}^N$ is much smaller
than $1/\epsilon$, the gap between the two first bands of the periodic
problem on the line, where
\beq\label{defmB} m_{B,\Omega}^N  = \inf_{\sum_j||\psi_{j,\perp}||^2=1} \Epsilon_{B,\Omega}^N
((\psi_{j,\perp}))\,. \eeq
An independent difficulty  is then to have more accurate estimates $m_A^N$ and $m_{B,\Omega}^N$
according to the regime of parameters. We do not have universal
estimates for this but have to separate two cases~:\begin{itemize}
\item the Weak Interaction case, where the interaction term ($L^4$
term) is at most of the same order as the ground state of the linear
problem in the same direction;
\item the Thomas Fermi case, where the kinetic energy term is much
smaller than the potential and interaction terms.
\end{itemize}
In what follows, when $N$ is not mentioned in $m_A^N$, $m_{B,\Omega}^N$,
 $\Epsilon_A^N$, $\Epsilon_{B,\Omega}^N$, then the notations are for
 $N=1$.
 Similarly, if $\Omega$ is not mentioned, this means that either the
 considered quantity is independent of $\Omega$ or that we are
 treating the case $\Omega=0$. To mention the dependence on other
 parameters, we will sometimes explictly write this dependence
 like for example $m_{A}^N ( \epsilon,\widehat g)$ or
 $m_{B,\Omega}^N(\tilde g, \omp)\,$.\\

\subsubsection{Case ($A$) : the longitudinal model} We consider  states which are
of type (\ref{tfa})
  with $\varphi \in L^2(\mathbb R_z/(NT)\mathbb Z)$.
 The energy of such test functions provides the upper bound
\beq\label{easyubA} E_\Omega ^{per,N} \leq \omega_\perp + m_A^N(\epsilon,
\widehat g)\; \eeq where $m_A^N$ is given by (\ref{defmA}) and
$\widehat g$
 was introduced in \eqref{defwidehatg}.

In order to show that the upper bound is an approximate
 lower bound, we first address the ``Weak
Interaction'' case,
 \beq
\label{AsscasAa} (AWIa)\quad 1<< \epsilon (\omega_\perp -\Omega)\,,\eeq
 and,  for a given $c>0$,
\beq\label{AsscasAb}
 (AWIb)\quad
g\omega_\perp  \epsilon^\frac 12 \leq c\,.
\end{equation}
The first assumption implies that  the lowest eigenvalue
$\lambda_{1,z}$
 of  the linear problem
in the $z$ direction (having in mind \eqref{asymp1})
 is much smaller than the gap in the transverse direction
 $\delta_\perp= \omp-\Omega$. This will allow the projection
onto the subspace  $\psi_\perp \otimes L^2(\mathbb R_z/(NT)\mathbb
Z)$. The second assumption implies that
 the nonlinear term (of order $g\omp/\sqrt\epsilon$) is of the same order as $\lambda_{1,z}$.
 It implies using \eqref{asymp1}, \eqref{asymptIN}  and the universal estimate
\begin{equation}
\lambda_{1,z} \leq m_A^N \leq \lambda_{1,z} + I_N\,,
\end{equation}
 that \beq \label{approxma}
 m_A^N\approx \frac 1 \epsilon\,.
\eeq
Here $\approx$ means ``of the same order'' in the considered regime of
 parameters. More precisely we mean
by writing \eqref{approxma} that, for any $\ep_0>0$,  there exists $C>0$
such that, for all $\ep \in ]0,\ep_0]$, any $g, \omp$ satisfying
 \eqref{AsscasAb},
$$
 \frac {1}{C \ep} \leq  m_A^N \leq \frac {C}{\ep}\,.
$$
Note that most of the time, we will not control the constant with
respect to $N$.\\
 All these rough estimates
are obtained by rather elementary semi-classical methods which are
recalled in Section \ref{Section1D}. More precise asymptotics of $
m_A^N$ will be given under the
 additional Assumption \eqref{condbiza} in Section
\ref{sssB12}.
 Thus, by
\eqref{AsscasAa}, $m_A^N$  is  much smaller than $\delta_\perp$. We will
prove
\begin{theorem}\label{theoA}~\\
When $\epsilon$ tends to $0$, and under Conditions (\ref{AsscasAa})
and (\ref{AsscasAb}), we have
\begin{equation}
E_\Omega^{per,N}=\omega_\perp+ m_A^N(\epsilon, \widehat g ) \; (1+o(1))\,.
\eeq
\end{theorem}
~\\
We now describe the  ``Thomas-Fermi'' regime, where we can also
justify the reduction to the longitudinal model.
 We assume that, for some given  $c>0\,$,
 \beq
\label{ATFa} (ATFa)\quad  g \omp \sqrt{\ep} >>1\,, \eeq

\beq \label{ATFb} (ATFb)\quad  g \omp \ep^2 \leq c \,, \eeq

\beq \label{ATFc} (ATFc)\quad  g^\frac {5}{12} \ep^{-\frac
16}\omp^{\frac {5}{12}} << (\omp -\Omega)^{\frac 38}\,. \eeq Note
 that \eqref{ATFa}  is the
converse of (\ref{AsscasAb}) while \eqref{ATFa} and \eqref{ATFc}
 imply that $1<<\ep(\omp-\Omega)$. This  implies 
 $\lambda_{1,z}<<\delta_\perp$, which is  the main condition to reduce to case
 A.
 Assumptions \eqref{ATFa} and
\eqref{ATFb} allow
 to show that~:
\beq \label{approxmATF}
  m_A^N \approx   \left(\frac{g\omp}{\ep}\right)^\frac 23\,,
\eeq
and this also implies
 that the nonlinear term is much bigger than $\delta_z$. \\
The estimate \eqref{approxmATF} will be shown in Section
\ref{ssTFr}, together with more precise ones with stronger
hypotheses (see
 Assumption \eqref{assthfb} and  \eqref{assmA}).

\begin{theorem}\label{theoATF}~\\
When $\epsilon$ tends to $0$, and under Conditions  \eqref{ATFa},
\eqref{ATFb} and \eqref{ATFc}, we have, as $\ep \ar 0$,
\begin{equation}
E_\Omega^{per,N}=\omega_\perp+ m_A^N(\epsilon, \widehat g ) \; (1+o(1))\,. \eeq
\end{theorem}
The proofs give actually  much stronger results.

\subsubsection{Case ($B$) : the transverse model}
 This corresponds to
 the idea of a  reduction on the ground eigenspace  in the $z$ variable,
 where the interaction term
 is kept in the transverse problem: therefore, this is a regime where
  $\omp\ep<<1$. We recall that we denote by
 $\lambda_{1,z}$ the ($N$-independent) ground state energy of $H_z^{per,N}$ and by
 $\phi_1^N$ the normalized ground state.
We consider  states which are of type (\ref{tfb}) or  \eqref{tfbN}.
 We have defined $\Epsilon_{B,\Omega}^N$ by
  \eqref{defEBN}-\eqref{defPsiN} and
 $m_{B,\Omega}^N$,
 the infimum of the energy of such test functions by \eqref{defmB}.
We have the upper bound \beq\label{univubB} E_{\Omega}^{per,N} \leq
\lambda_{1,z}+  m_{B,\Omega}^N\,. \eeq When $N=1$, $m_{B,\Omega}$ is a
function of $\tilde g$ and $\omp$
 as it is clear from \eqref{EB}
 and \eqref{defmB}.
 Note that, from \eqref{roughphi4}, we get
\begin{equation}
\tilde g=  g  (\int_{-\frac T2}^{\frac
T2} \phi_1(z)^4 dz)\approx \frac{g}{\sqrt\ep}\,.
\end{equation}
Again we can discuss two different cases according to the size of
the interaction. In the Weak Interaction case, we prove the
following~:
\begin{theorem}\label{theoBWI}~\\
When $\epsilon$ tends to $0$, and under the conditions
\begin{equation}\label{BWIa}
(BWIa)\quad g \epsilon^{-\frac 12}  \leq C \,, \eeq \beq
\label{BWIb}(BWIb) \quad \omp \epsilon  << 1\,,
\end{equation}
 then \beq E_{\Omega}^{per,N}=\lambda_{1,z}+ m_{B,\Omega}^N
 (1+o(1))\,.
\eeq
\end{theorem}
Condition (BWIb) implies that the bottom of the spectrum of the
linear problem in the $x-y$ direction is much smaller than
$\delta_z$, the gap in the $z$ direction, which is of order $1/\ep$.
Condition \eqref{BWIa}, together with \eqref{universal} and \eqref{asymptIN}, implies that
 $m_{B,\Omega}^N$ satisfies
\beq m_{B,\Omega}^N\approx \omp\,. \eeq Indeed, (BWIa) and (BWIb)
 imply $g \epsilon^{\frac 12}  \omp  <<1$, that is (QLb).

In the Thomas-Fermi case, we prove the following~:
\begin{theorem}\label{theoBTF}~\\
When $\epsilon$ tends to $0$, and under the conditions
\begin{equation}\label{asscasBa}
(BTFa)\quad \sqrt \epsilon  << g \,, \eeq \beq
\label{asscasBb}(BTFb)\quad \omp \sqrt{g} \epsilon^\frac 34  << 1\,,
\end{equation}
and\beq\label{asscasBc} (BTFc)\quad g^\frac 32 \epsilon^\frac 14
\omp   << 1\,,\eeq then \beq E_{\Omega}^{per,N}=\lambda_{1,z}+ m_{B,\Omega}^N
 (1+o(1))\,.
\eeq
\end{theorem} Note that $(BTFa)$ is the converse of $(BWIa)$.
We will see in Proposition \ref{estmB} (together with \eqref{majcasB},
 \eqref{asympcasBOmegaN=1}
and  \eqref{lbmbntf}) that, under these
assumptions and Assumption \eqref{controleOmega},
 the term $m_{B,\Omega}^N$ satisfies
\beq
m_{B,\Omega}^N\approx \omp\sqrt g/\ep^{1/4}\,,
\eeq
 and
thus is much smaller than  $\delta_z^N$ which is  of order $\frac 1
\epsilon$.

Our proofs are made up of two parts~: rough or accurate  estimates of $m_{A,\Omega}^N$
and $m_{B,\Omega}^N$ on the one hand
 and a lower bound for $E_\Omega^{per,N}$ on the other hand. The
lower bound consists in showing that the upperbound obtained by  projecting on the special states introduced
above in \eqref{tfa}, \eqref{tfb} or \eqref{tfbN} is actually also
 asymptotically a good lower bound.

\subsubsection{Tunneling effect and discrete model}

 Since the Wannier functions are localized in the $z$ variable,
 the energy of a function
 $\Psi=\sum_{j=0}^{N-1}\psi_j^N(z)\psi_{j,\perp}(x,y)$ provides at leading order
 the sum of $N$ decoupled energies for $\psi_{j,\perp}$ on each slice $j$.
 At the next order,
  in the computation of the $L^2$ norm of the gradient, only the nearest
 neighbors in $z$ interact through an exponentially small term, describing
 what is called the tunneling effect.
  These simplifications are discussed in section \ref{stunneling}.
  We are lead to new functionals and in particular a discrete model
  that we analyze in relationship with the physics papers.

  In case A, the behavior on each slice $j$ is the same, given by
  $\psi_\perp$ and it is the behavior on the $z$ direction which
  has a tunneling contribution. There are no vortices whatever the
  velocity $\Omega$.

  In case B, for $N=1$, there are vortices for large velocity
  and they are located on each slice at the same place.
   For $N$ large, it is an open and interesting question to analyze
   whether it is possible for a vortex line to vary location
   according to the slice, whether vortices
   interact between the slices and how. This could be performed using our reduced
   models.

\subsection{Organization of the paper}
The paper is organized as follows.
 In Section \ref{s2}, we start the spectral analysis of the
linear problems in  the longitudinal  and  transverse
 directions. We recall in particular the main techniques which can be
 used for the analysis of the spectral problem with periodic potential
 on the line. Section \ref{Section1D} is devoted to the semi-classical results
 for the periodic problem. Although we are mainly interested in
 $1D$-problems
 we recall here techniques which are true in any dimension and can be
 useful for the analysis of $2D$ or $3D$ optical lattices.

In Section \ref{scasA}, we prove  the main theorems for case A. In
Section \ref{s1D}, we analyze the ground
state of the $1D$ nonlinear
 energy $\Epsilon_A^N$ for $N=1$ and $N>1$ and also
 distinguish between the two cases: Weak Interaction and Thomas-Fermi.
Section \ref{scasB} corresponds to a similar analysis for the
transverse models $\Epsilon_B^N$. Section \ref{stunneling} is
devoted to the tunneling effects and discuss, on the basis of the
semi-classical estimates
 of Section \ref{Section1D}, some
 results by physicists on the discrete nonlinear Schr\"odinger model.
  In Section \ref{s4}, we  analyze various boundary conditions and
   compare in particular the problems on $\mathbb R^3$ (which is
 completely solved) and the problems on
 $\mathbb R^2 \times ]-\frac{NT}{2}, \frac{NT}{2}[$ with periodic condition which seem
 physically more interesting.

\section{Analysis of the linear model}\label{s2}
The linear model which appears naturally is associated to
$$
H_\Omega = H^\Omega_{\perp} + H_z\,,
$$
which was presented in the introduction (see
\eqref{3.2}-\eqref{3.6}).
  A natural
condition (for the strict positivity of the operator
$H_{\perp}^\Omega$) is
 Condition \eqref{conditionsurOmega}.
In this situation with separate variables, we can split the spectral analysis
 in the separate spectral
analysis of  $H_{\perp}^\Omega$ and the spectral analysis of a
suitable realization of $H_z$ which will be
 presented in the next subsection.

\subsection{The harmonic oscillator in the transverse variable}

For simplicity,
 we begin the analysis of $H_\perp=H_{\perp}^\Omega$ with  the case
$
\Omega =0$.\\
 The first
operator $H_{\perp}^{\Omega=0}$ is a harmonic oscillator with
discrete spectrum and the bottom of its spectrum
 is given by
\begin{equation}
\inf (\sigma(H^0_{\perp} )) =   \omp \,.
\end{equation}
A corresponding $L^2$-normalized ground state is the Gaussian
\begin{equation}
\psi_{\perp} = \left(\frac{\omp}{\pi}\right) ^\frac 12 \exp - \frac {\omp}{2} r^2\,.
\end{equation}
Moreover the gap between the ground state energy and the second
eigenvalue (which has multiplicity $2$) is given by
\begin{equation}
\lambda_{2,\perp} - \omp =  \omp \,.
\end{equation}

The spectrum of $H_\perp^\Omega$  can be recovered by considering first
the joint spectrum of $H^\Omega_\perp$ and  $L_z$. For each eigenspace
of $L_z$ corresponding to $ \ell$ for some $\ell \in \mathbb Z$,
 we can look at the operator
\begin{equation}
H_\perp^{(\ell)}:=- \frac{1}{2} \Delta_{x,y} +
 \frac 12 \omp^2 r^2 -
 \Omega \ell  \,,
\end{equation}
considered as an unbounded operator
 on $L^2(\mathbb R^2)\cap \mbox{ Ker }(L_z -  \ell)$.

More precisely, for $\Omega$ satisfying \eqref{conditionsurOmega}, a
common eigenbasis of $L_z$ and $H_\perp^0$ is given by the set of
(not normalized) Hermite functions:
  \begin{equation}
\phi_{j,k}( x,y)=e^{\frac{\omp}{2}(x^2+y^2)}\;
(\partial_x+i\partial_y)^j\;(\partial_x-i\partial_y)^k
\left(e^{-\omp (x^2+y^2)}\right)
 \label{Hermite}
  \end{equation}
where $j$ and $k$ are non-negative integers.\\
 The eigenvalues are
$(j-k)$ for $L_z$ and
 \begin{equation}
 E_{j,k}=\omp+ (\omp-\Omega)j+ (\omp+\Omega)k
\label{spectrum}
 \end{equation}
for $H_\perp^\Omega$.

The spectrum of $H_\perp^{(\ell)}$ is obtained by considering the pairs
$(j,k)$ such that $j-k=\ell$.\\
We emphasize that this orthogonal  basis of eigenfunctions
 is independent of $\Omega$.

\subsection{The band spectrum in the longitudinal direction}
 The
second operator $H_z$ can be analyzed by semi-classical methods but
note that our semi-classical parameter is $\ep$. One can of
course in the case of the specific $w$ introduced in \eqref{1.5} recognize this operator
 as the Mathieu operator (for which a lot of information
 can be obtained using special functions (see \cite{AS}))
 but we prefer to give the presentation of
the theory for a more general periodic
 potential
  $w$. We hope that the general ideas which are behind will become clearer.

There are two related approaches for the analysis of the spectrum of
$H_z$, which is known to be a band spectrum, i.e.  an absolutely
continuous spectrum which is  a union of closed intervals, which are
called the bands.
\subsubsection{Floquet's theory}
We can first  use the Floquet theory (or the Bloch theory, which is
an
 alternative name for the same theory). This is more detailed in the
 appendix.  One can show that
 the spectrum of $H_z$
 is obtained by taking the closure
 of $\cup_{k\in [0,2\pi/T]}\sigma( H_{z,k})$
 where
$$
H_{z,k}=   - \frac{1}{2} \left(\frac{d}{dz} + ik\right)^2 + W_\ep(z)
$$
is considered as an operator on $L^2(\mathbb R/T\mathbb Z)$. So
\begin{equation}
\sigma (H_z) =\overline{\cup_{k\in [0,\frac{2 \pi}{T}]}
\sigma(H_{z,k})}\,.
\end{equation}We now
write
\begin{equation}
\Gamma = T\mathbb Z \mbox{ and } \Gamma^* =  \frac{2\pi}{T}\mathbb
 Z\,.
\end{equation}
Hence we have to
analyze for each $k$ the operator $H_{z,k}$ on  $L^2(\mathbb R/\Gamma)$.
Later we will use the notation
\beq \label{defHzper}
H_z^{per}=H_{z,0}\,.
\eeq

A unitary equivalent presentation of this approach consists in
analyzing $H_z$
 restricted to the subspace $\mathfrak h_k$  of the $u\in L^2_{loc}(\mathbb R)$
 such that
\begin{equation}\label{Fl}
u(z+T) = e^{ ik T  }\, u(z)\,.
\end{equation}
Here  we did not see a $k$-dependence in the differential operator
 but this is the choice of the space $\mathfrak h_k$ (which is NOT in
 $L^2(\mathbb R)$), which gives the
 $k$-dependence. Condition \eqref{Fl} is called a Floquet
 condition.\\
This means that we have written, using the language of the
 Hilbertian-integrals, the decomposition
\begin{equation}\label{HilDec1}
L^2(\mathbb R) =\int^{\oplus}_{[0,2\pi/T]}
 \mathfrak h_k    \; dk
\end{equation}
 and  that we have for the operator the corresponding decomposition
\begin{equation}\label{HilDec2}
H_z =\int^{\oplus}_{[0,2\pi/T]} \widetilde H_{z,k}\; dk \,,
\end{equation}
with $\widetilde H_{z,k}$ unitary equivalent to $H_{z,k}$.\\

For each $k\in [0,2\pi/T[$,   $H_{z,k}$ has a discrete spectrum which can be described
 by an increasing
 sequence of eigenvalues $(\lambda_{j}(k))_{j\in \mathbb N}$. The spectrum of $H_z$
 is then a union of bands $B_j$, each band being described by
 the range of $\lambda_j$. At least when we have the additional symmetry
 $W_\ep$ even, one can determine for which value of $k$ the
 ends of the band $B_j$ are  obtained. For $j=1$, we know in addition from the
 diamagnetic inequality  that the
  minimum of $\lambda_1$ is obtained for $k=0$~:
\beq
\inf_k \lambda_1(k)=\lambda_1(0)\,.
\eeq

\subsubsection{Wannier's approach}
When the band is simple (and this will be the case for the lowest band
 in the regime $\ep$ small), one can associate to $\lambda_j(k)$
 a normalized\footnote{in $L^2(]-\frac T2,\frac T2 [)$,}   eigenfunction $\varphi_j(z,k)$ with in addition an analyticity
 with respect to $k$ together with the $(2\pi/T)$-periodicity in $k$.

In this case (we now take $j=1$), one can associate
 to $\varphi_1$, which satisfies,
\begin{equation}\label{Wa0}
\varphi_1(z+T;k) =\varphi_1(z,k)\,,
\end{equation}
and
\begin{equation}\label{Wa1}
\varphi_1(z;k+\frac{2\pi}{T}) =\varphi_1(z,k)\,,
\end{equation}
 a family of  Wannier's functions
$(\psi_\ell)_{\ell \in \Gamma}$ defined by
\begin{equation}\label{defpsi0}
\psi_0(z) = \frac{T}{2\pi}\int_0^{\frac{2\pi}{T}} \exp   (i kz)\;
\varphi_1(z,k)\,dk\,,\; \psi_\ell(z)=\psi_0(z-\ell
)\,,\end{equation}
 for $\ell \in \Gamma $\,.\\
In addition, we can take  $\psi_0$  real. One can indeed
 construct  $\varphi_1$ satisfying
 in addition the condition
\begin{equation}\label{flo1}
\overline{\varphi_1(z,k)} = \varphi_1(z, -k)\,.
\end{equation}
One obtains  (after some normalization of $\psi_0$) that
\begin{proposition}~
\begin{enumerate}
\item
 The family $(\psi_\ell)_{\ell\in \Gamma}$ gives  an orthonormal basis of the spectral space attached to the first
 band.
\item
  $\psi_0$ is an exponentially decreasing function.
\end{enumerate}
\end{proposition}
 The second point can be proved using  the
 analyticity\footnote{One can make a contour deformation in the
 integral defining  $\psi_0$ in \eqref{defpsi0}.} with respect to $k$.\\
This orthonormal basis corresponding to  the first band
 plays  the role of the basis  $P_j(z)\exp -\frac{|z^2|}{2}$
 in the Lowest Landau Level approximation.\\
Note that we recover $\varphi_1(z,k)$
 by the formula
\begin{equation}\label{flo2}
\varphi_1(z,k)= \exp (- i kz) \;\sum_{\ell \in \Gamma} \exp  (ik\ell
) \;\psi_\ell (z)\,.
\end{equation}
Moreover, the operator $A$ on $\ell^2(\Gamma )$
 whose matrix is given by
\begin{equation}\label{flo3}
A_{\ell\ell'} = \langle H_z \psi_\ell,\psi_{\ell'}\rangle
\end{equation}
is unitary equivalent to the restriction
 of $H_z$ to the spectral space attached to the first band.\\
One can of course observe that $A$ commutes with the translation
 on $\ell^2(\Gamma)$, so it is a convolution operator by a sequence
 $a\in \ell^1(\Gamma)$ (actually in the space of the rapidly
 decreasing sequences ${\small {\mathcal S}}(\Gamma)$),
\begin{equation}\label{flo4}
A_{\ell\ell'}=a(\ell -\ell')\,,
\end{equation}
 which is actually
 the Fourier series of $k\mapsto \lambda_1(k)$
\begin{equation}\label{flo5}
\widehat \lambda_1=a\,,
\end{equation}
where
\begin{equation}\label{flo6}
\widehat \lambda_1 (\ell):= \frac {T}{2\pi} \int_0^{2\pi/T} \exp (-i
\ell k)\; \lambda_1(k)\; dk\,.
\end{equation}
So we have
$$
(Au)(\ell)=\sum_{\ell'\in \Gamma} a(\ell -\ell')
  u(\ell')\,,\mbox{ for } u\in \ell^2(\Gamma)\,.
$$
\subsubsection{($NT $)-periodic problem}\label{sssNper}
There is another way to proceed at least heuristically. We keep $w$
$T$-periodic
 but look at the
$(NT)$-periodic problem  and we analyze this problem. The
spectrum is discrete but the idea is that we will recover the band
spectrum
 in the limit $N\rightarrow +\infty$. If we compare with what we do in
 the Floquet theory, the analysis of the $(NT)$-periodic problem
 consists in considering the direct sum of the problems with a Floquet
 condition corresponding
 to $k =0, \frac{2\pi}{N T },\cdots, \frac{2\pi (N-1)}{N T}$.\\

Note that this decomposition into a direct sum works only for linear
problems,
 so it will be interesting to explore this approach for the non linear
 problem.

In this spirit, it can be useful to have an adapted orthonormal basis
 of the spectral space attached to the first $N$ eigenvalues of the
 $NT$-periodic
 problem (which can be identified with the vector space generated
 by
 the eigenfunctions corresponding to the $N$ Floquet eigenvalues
 associated with
  $k =0, \frac{2\pi}{NT},\cdots, \frac{2\pi (N-1)}{NT }$.\\

Our claim is that there exists an orthonormal  basis, for the
$L^2$-norm on $]-\frac{NT}{2},\frac{NT}{2}[$,  consisting of $(NT)$-periodic
functions and  replacing the Wannier functions.\\

We write
\begin{equation}\label{Nper1}
\psi^N_0(z) = \frac{1}{\sqrt{N}} \sum_{j=1}^{N} \phi^N_j(z)\,,
\end{equation}
where $\phi^N_j$ is an eigenfunction\footnote{Note that except in
the
 case $j=1$, we do not claim that $\phi^N_j$ is the $j$-th
 eigenfunction
 but this is the first one corresponding to the condition \eqref{Nper2}.}
 of the $(NT)$-periodic problem, chosen in such a way that
\begin{equation}\label{Nper2}
\phi^N_j(z+T)= \omega_N^{j-1} \phi^N_j(z)\,,
\end{equation}
with $\omega_N = \exp (2 i \pi/N)$.\\
We can then introduce
\beq
\Gamma^N=\Gamma/(NT\mathbb Z)\,,
\eeq
and  define, for $\ell \in \Gamma^N=\Gamma$, the $(NT)$-Wannier
functions
\begin{equation}\label{Nper3}
\psi^N_\ell (z) =\psi_0^N(z-\ell)\;
\end{equation}
This gives an orthonormal basis of the eigenspace attached
 to the first $N$  eigenvalues of the $(NT)$-periodic problem.
 These first $N$ eigenvalues belong to the previously defined first
 band.\\

Note that conversely, we can recover the eigenfunctions $\phi_j^N$
 from the $\psi_j^N$ by a discrete Fourier transfrom. In particular we
 have
\beq\label{deWaversphi}
\phi_1^N = \frac{1}{\sqrt{N}} \sum_{j=0}^{N-1} \psi_j^N\,.
\eeq

Except the fact that these ``Wannier'' functions
are NOT exponentially decreasing at $\infty$ (they are by construction
$(NT)$-periodic), one can then play with them in the
same way (this corresponds to the replacement of the Fourier
 series by the finite dimensional one). We then meet
the ``discrete convolution'' on $\ell^2(\Gamma^N)$~:
$$
(A^Nu)(\ell) = \sum_{\ell'\in \Gamma^N} a_N(\ell -\ell')
u(\ell')\,,\mbox{ for }
 u\in \ell^2(\Gamma^N)\,.
$$
Of course $\ell^2(\Gamma^N)$ is nothing else than $\mathbb C^N$
 with its natural Hermitian structure.

We have presented different techniques to determine the bottom
 of the spectrum of $H_z$, which all provide the same ground energy.
We will now  recall more quantitative results
 based on the so-called semi-classical analysis.

\section{Semi-classical analysis for the $T$-periodic case}\label{Section1D}
\subsection{Preliminary discussion}
Till now, we have not strongly used that we are
 in a semi-classical regime: our semi-classical parameter
 here will not be the Planck constant $\hbar$ (which was
 already assumed to be equal to $1$) but $\ep$. We will now use
 this additional assumption for extracting quantitative results
 from the previously presented qualitative theory.
As already said, the physics literature is analyzing a very
particular model,
 the Mathieu equation. We will rapidly sketch how one can do this in
full generality. For the one dimensional case which is considered
here, one can probably refer to Harrell \cite{Ha} (who uses
techniques
 of ordinary differential equations) or to the book of Eastham \cite{Eas}, but we will describe a proof
 which is more general in spirit,
 which is not limited to the one dimensional situation (see Simon
 \cite{Si}, Helffer-Sj\"ostrand \cite{HeSj}, Outassourt \cite{Ou}) and is described in the books of Helffer \cite{He}  or
 Dimassi-Sj\"ostrand \cite{DiSj}.\\
As we have shown in the previous section, the
 good description of the first  band, can be either obtained by
 a good approximation of $\lambda_1(k)$ and $\varphi_1(z,k)$ as $\ep\rightarrow 0$
 or by first finding a good approximation of the Wannier function
 $\psi_0$ introduced in \eqref{defpsi0}, which is expected to be
 exponentially localized
 in one well, or of the $(NT)$-periodic  Wannier function
 introduced in \eqref{Nper1}.\\

The analysis is done usually in two steps. First
we localize roughly $\lambda_1(k)$, then
 we analyze very  accurately the variation of  $\lambda_1(k)-\lambda_1(0)$.\\
The first one will be obtained by a harmonic approximation
 and the second one by the analysis of the tunneling effect.\\

\subsection{The harmonic approximation}\label{ssha}~\\
We will provide the explanation in  a general case containing  the
model considered  by Snoek \cite{Sn} as a particular case. We recall
that we work under Assumption \ref{assw}.
 The statements below  are sometimes  written vaguely and we refer to
 \cite{DiSj}
 or \cite{He} for more precise mathematical statements.\\
For the approximation of $\lambda_{1,z}(0)$ (actually for any $\lambda_{1,z}(k)$)
 the rule is that we replace ${\bf w}(z)$ (having in mind \eqref{1.5}) by its quadratic approximation
 at $0$.
The harmonic  approximation consists in first looking at the operator
\begin{equation}\label{ha1}
- \frac{1}{2}\frac{d^2}{dz^2}  + \frac{{\bf w}''(0)}{2 \ep^2}z^2 \,,
\end{equation}
on $\mathbb R$.\\
For the model in \cite{Sn}, ${\bf w}(z)=\sin^2(\frac{\pi z}{T})$, and we
find
\begin{equation}\label{ha1sn}
- \frac{1}{2}\frac{d^2}{dz^2}  + \frac{1}{\ep^2} (\frac{\pi z}{T})^2\,.
\end{equation}
This operator is a harmonic oscillator whose spectrum is explicitly
known. The $j$-th eigenvalue is given by
\begin{equation}\label{eigenvHz}
\lambda_{j,z}^{har} = \frac{j-\frac 12}{\ep}\sqrt{{\bf w}''(0)}\,.
\end{equation}
The two main pieces of information we have to keep in mind are that
the ground state
 energy is
\begin{equation}\label{bottomHz}
\lambda_{1,z}^{har} = \frac{1}{2\ep}\sqrt{{\bf w}''(0)}\,,
\end{equation}
and that the gap between the first eigenvalue and the second value is
 given by
\begin{equation}\label{gapHz}
\delta^{har}_z:=\lambda_{2,z}^{har}-\lambda_{1,z}^{har}= \frac{1}{\ep}\sqrt{{\bf w}''(0)}\,.
\end{equation}
The corresponding positive $L^2$ normalized ground state
 is then given by
\begin{equation}\label{ha3}
\psi^{har} (z) = \pi^{-\frac 14} {\bf w}''(0)^{\frac 18}\,\ep^{-\frac 14}
\exp - {\bf w}''(0)^\frac 12  \frac{z^2}{2\ep}\,.
\end{equation}
It will also be important later to have the computation of the $L^4$
norm. So
 we get by immediate computation~:
\begin{equation}\label{normel4}
\int_\mathbb R \psi^{har} (z)^4\, dz =  \pi^{-\frac 12} {\bf w}''(0)^{\frac 14}\,\ep^{-\frac 12} \,.
\end{equation}
The mathematical result is that this value provides a good
approximation of $\lambda_{1,z}(0)$ (and hence of the bottom of the
spectrum of $H_z$) with an error which is $\Og(1)$ as $\ep
\rightarrow 0$~:
\begin{equation}\label{ha5}
\lambda_{1,z}(0) = \lambda_{1,z}^{har} + \Og (1)\,.
\end{equation}
By working a little more, one can actually obtain a complete
expansion of $\ep \lambda_{1,z}(0)$ in powers of $\ep$
 and hence, of $\ep
 \lambda_{1,z}(k)$, since they have the same expansion.
For each $j\in \mathbb N^*$, one has a similar expansion for
 $\ep \lambda_{j,z}(0)$. This implies in particular an estimate of
 $\lambda_{2,z}(0)-\lambda_{1,z} (0)$, called the longitudinal gap~:
\begin{equation}\label{longgap}
\delta_z :=\lambda_{2,z}(0)-\lambda_{1,z} (0)=
\frac{\sqrt{{\bf w}''(0)}}{\ep}+ \Og(1) \,\,.
\end{equation}\\ From now on,
we simply write $\lambda_{1,z}$  or $\lambda_1$ instead of $\lambda_{1,z}(0)$ for the
ground state energy of the periodic problem.\\
Let us note that the ground state of the harmonic oscillator
 also provides a good approximation of the ground state of $H_z^{per}$. So we obtain,
  using \eqref{normel4} that for $\phi_1$, the $L^2$-normalized ground state
of $H^{per}_z$, we have
\begin{equation}\label{normel4a}
\int_{-\frac T2}^{+\frac T2} \phi_1(z)^4 dz =  \pi^{-\frac 12}
{\bf w}''(0)^{\frac 14}\,\ep^{-\frac 12} + \Og(1)\,.
\end{equation}

\subsection{The tunneling effect}\label{sstunneling}~\\
We now briefly explain the results about the length of the first
band, which is exponentially small
 as $\ep \ar 0$.
The results can take the following form (see the work of Outassourt
 \cite{Ou}
 or the book by Dimassi-Sj\"ostrand, Formula (6.26))
\begin{equation}
\lambda_1(k) -\lambda_1(0) =  2 (1-\cos (k T) ) \tau + \Og (\exp -\frac{S+\alpha}{\ep})
\end{equation}
with $\alpha >0$ (arbitrarily close from below to $1$) and, for some $c_\tau \neq 0$,
\begin{equation}\label{asympoft}
\tau \sim  c_\tau \; \ep^{-\frac 32} \; \exp -\frac{S}{\ep}\,.
\end{equation}
Moreover one can express the constants $c_\tau$ and $S$  once ${\bf w}$ is given
(see\footnote{The computation
 is a little simpler in the case when ${\bf w}$ is
 even.} also \cite{He} in addition to the previous references).
This $\tau $ seems  to be called  in some physical
 literature
 the hopping amplitude.\\
Here, we simply explain how one computes $S$ which determines the
exponential decay of $\tau $ as $\ep \ar 0$. In any
 dimension, $S$ is interpreted
 as the minimal Agmon distance between two different minima of the
 potential $w$. In one dimension, with $w$ satisfying Assumption
 \eqref{assw},
 this distance is
 simply the Agmon distance between two consecutive minima and is given
 by
\begin{equation}
S:=\sqrt{2}\; \int_{-\frac T2}^{\frac T2} \sqrt{{\bf w}(z)} \, dz\,.
\end{equation}
In particular, when ${\bf w}(z) = \sin^2(\frac{\pi z}{T})$, we get
\begin{equation}
S:= \sqrt{2}\;\int_{-\frac T2}^{\frac T2}| \sin (\frac{\pi
  z}{T})|\, dz =  \frac{2 \sqrt{2}T}{ \pi}\,.
\end{equation}
This is to compare to (14) in \cite{SnSt}, which is not an exact
formula (as wrongly claimed) but only  an asymptotically correct
formula. It  can be found, for this Mathieu operator, in
\cite{AS}.\\
Let us give the formula for the constant $c_\tau$. It can be found
in \cite{Ha}, see also \cite{Ou}, Formula (4.14) and \cite{He}
p.~58-59. We have~:
\begin{equation}\label{ctau}
c_\tau = 2^{\frac 34} \pi^{-\frac 12} \exp  A_\tau \,,
\end{equation}
with (assuming $w$ even)
\begin{equation}
A_\tau  =  \lim_{\eta \ar 0} \left( \int_\eta^{\frac T2} \frac{1}{
\sqrt{{\bf w}(z)}}\, dz+ \frac{\sqrt{2}}{\sqrt{{\bf w}''(0)}} \ln
  \eta \right)\,.
\end{equation}
We just sketch the mathematical proof. Filling out all the wells suitably
 except one (say $0$), we get a new potential $w^{mod}\geq {\bf w}$ which
 coincides
 with ${\bf w}$ in an interval containing $0$ and excluding small
 neighborhoods of all the other  minima. We consider,
 for $\ep$ small enough, the ground state of this modified problem and
 (multiplying by a cut-off function)
 we get a function $\psi_0^{app}$ (and an eigenvalue $\lambda_1^{app}$) which is a very good approximation
 of  $\psi_0$.\\
Now the hopping amplitude in the abstract theory is given\footnote{For
 the Mathieu potential, this  is consistent with
Formula (13) in \cite{SnSt}.} {\bf exactly} by
\begin{equation}
 - \tau =a(T) =  \langle H_z \psi_0 \,,\; \psi_1\rangle =  \langle (H_z -\mu)\psi_0 \,,\; \psi_1\rangle\,,
\end{equation}
the last equality being satisfied, due to the orthogonality of
$\psi_0$ and $\psi_1$, for any $\mu$.   When replacing $\psi_0$  by
its approximation, one has to be careful, because $\psi^{app}_0$ and
$\psi^{app}_1:=\psi^{app}_0 (\cdot -T)$ are no more orthogonal. So
this leads to  take $\mu = \lambda_1^{app}$, and one can prove that
\begin{equation}\label{approxt}
\tau  \sim - \langle (H_z -\lambda_1^{app})\psi^{app}_0 \,,\; \psi^{app}_1\rangle\,.
\end{equation}
An easy way to see that $\tau$ is exponentially small is to observe
that
\begin{equation} \langle (H_z -\lambda_1^{app})\psi^{app}_0 \,,\;
\psi^{app}_1\rangle
 = \ep^{-2}\;  \langle ({\bf w}(z) - w^{mod} ) \psi^{app}_0 \,,\;
 \psi^{app}_1\rangle\,,
\end{equation}
 and to use the information on the asymptotic decay of $\psi^{app}_0$.
The WKB-approximation of $\psi^{app}_0$ is,  in a neighborhood of
$0$,  \beq \label{approxwkb} \psi^{wkb}_0= \ep^{-\frac 14} \,
b(z,\ep) \exp -  \frac 1 \ep \int_0^z\sqrt{{\bf w}(s)} ds \,,\, \mbox{ for
} z\geq 0\,, \eeq with \beq b(z,\ep) \sim \sum_{j\geq 0} b_j(z)
\ep^j\,, \eeq and \beq b_0(z)=\pi^{-\frac 14} \exp \left(-
\int_{0}^{z}\,
 \frac{(w^{\frac{1}{2}})'(t) - \sqrt{\frac{{\bf w}''(0)}{2}}}{2 \sqrt{{\bf w}(t)}}\,dt \right)  \,.
\eeq It should then be completed by symmetry to get an even WKB
solution on
 $]-T,+ T[$.\\
Note that we have
$$
(w^\frac 12)'(T_-) = -  \sqrt{\frac{{\bf w}''(0)}{2}}\,,
$$
which implies that $b_0$ tends to $+ \infty$ as $z\ar T_-$.\\

An integration by parts together with a WKB approximation  leads to
the
asymptotic  estimate of $\tau $ announced in (\ref{asympoft}).
More precisely, we get that the prefactor $c_\tau$ is immediately
related to the constant  $ b_0(\frac T2)^2\sqrt{V(\frac T2)}$ and this leads to
\eqref{ctau}. Note that more generally we have
\begin{equation}\label{constance}
b_0(z) b_0(T-z) \sqrt{V(z)} = \mbox{ Cst}\,,
\end{equation}
which again shows the blowing up of $b_0$ at $T$.\\

Finally, we emphasize that $\psi_0^{wkb}$ is  a good approximation
 of $\psi_0 $ only in intervals $]-T + \eta, T-\eta[$
 for some $\eta >0$.\\

One can also see that $a(kT)$ is of the order of $|a(T)|^{|k|}$ (for
$k\geq 2$)
\begin{equation}
a(kT) = \widetilde \Og (\tau^2)\,,
\end{equation}
so it is legitimate  in order to compute the
 width of the first band to forget all the $a(\ell)$ for $\ell\in
 \Gamma, \ell \neq 0,\pm T$.\\
Thus, in the $k$ variable, the spectrum (corresponding to the first
band) is up to a very small error, of the order of the square of
$a(T)$,
 given by the operator of multiplication in $L^2(\mathbb R/\Gamma)$ by
  the function $a(0) + 2 a(T) \cos( kT)$.\\
\begin{remark}~\\
What is written above corresponds to the use of Wannier functions on
 $\mathbb R$. One can write a close theory using the $(NT)$-periodic
 Wannier functions without  modifying the main terms of the
 asymptotics. In particular, $\psi_0^{wkb}$ is also  a good approximation
 of $\psi_0^{N}$   in intervals $]-T + \eta, T-\eta[$
 for some $\eta >0$.\\
The interest of the Wannier functions on $\mathbb R$ is that
 they allow to recover the information for all
 Floquet eigenvalues (see the discussion in Section \ref{sectionDNLS}).
\end{remark}




\section{Justification of the reduction to the longitudinal  energy
  $\Epsilon_A^N$}\label{scasA}
\subsection{Main result}
In this section, we address the reduction to the energy
$\Epsilon_A^N$ defined in \eqref{EpA} and  prove the following
theorem (recall that $m_A^N$ is defined in \eqref{defmA}):
\begin{theorem}\label{theoAfull} If
\beq\label{newAa} (A\Omega a) \quad \; m_A^N(\ep,\widehat
g)(\omp-\Omega)^{-1} << 1 \eeq
 and
\beq \label{newAb} (A\Omega b) \quad \; g (2\omp -\Omega)
m_A^N(\ep,\widehat g) (\omp-\Omega)^{-\frac 32}  << 1, \eeq
 we have
\begin{equation}\label{equivA}
\inf_{||\Psi||=1} \Epsilon_{\Omega}^{per,N}(\Psi) = \omp+
m_A^N(\ep,\widehat g)  (1+o(1)).
\end{equation}
\end{theorem}
Both Theorem \ref{theoA} and Theorem \ref{theoATF} are a consequence
of Theorem \ref{theoAfull} as soon as we have the appropriate
rough  estimates on $m_A^N$ already presented in the introduction. This
is what we explain first in Subsection \ref{subest} before proving
the theorem in Subsection \ref{subpr}.

\subsection{Proof of Theorem \ref{theoA} and Theorem
\ref{theoATF}}\label{subest}
\subsubsection{Weak Interaction case}\label{ssswima}
In the Weak Interaction case, we recall from \eqref{approxma},
 that, when \eqref{AsscasAb} is
satisfied, then
\begin{equation}\label{approxmAN}
m_A^N\approx  1 /\epsilon\,.
\end{equation}
 Therefore,
 when  \eqref{AsscasAa} and \eqref{AsscasAb} are satisfied, then
\eqref{newAa} and \eqref{newAb}
  automatically hold  with  the observation that
$$
g (2\omp-\Omega)(\omp-\Omega)^{-\frac 32} m_A^N(\ep,\widehat g)\leq
C
 g (2\omp-\Omega)\epsilon^{\frac 12}((\omp-\Omega)\epsilon)^{-\frac
 32}<<1\,,
$$
and   Theorem \ref{theoA} follows from Theorem \ref{theoAfull}.
\subsubsection{Thomas-Fermi case}\label{ssstfma}
In the Thomas-Fermi case,
 we will prove  in \eqref{rhf1} that,
 when  \eqref{ATFa} and \eqref{ATFb} are
satisfied, then
\begin{equation}\label{majTFa}
m_A^N\approx ( g \omp /
\epsilon)^{2/3}\,.
\end{equation}

Let us verify that, if  \eqref{ATFa},  \eqref{ATFb} and \eqref{ATFc} are satisfied, then
\eqref{newAa} and \eqref{newAb}
hold. This will prove Theorem \ref{theoATF}.

We get \eqref{newAa} in the following way.
First we have~:
$$
(\omp-\Omega)^{-1} m_A^N(\epsilon,\widehat g) \leq C (\omp-\Omega)^{-1}
 \omp^{\frac 23 }g^\frac 23
\epsilon^{-\frac 23}\,.$$ Hence \eqref{newAa} is a consequence of
\begin{equation}\label{ATFd}
g \omp << \epsilon (\omp-\Omega)^{\frac 32}\,,
\end{equation} which follows from \eqref{ATFc} since \eqref{ATFa} and \eqref{ATFc} imply that
$(\omp-\Omega) \epsilon >>1$. The check of \eqref{newAb} is then
immediate from \eqref{ATFc}
and \eqref{majTFa}.\\

\subsection{Proof of Theorem \ref{theoAfull}}\label{subpr}
Because of the upper bound (\ref{easyubA}), Theorem \ref{theoAfull}
is a consequence of the following proposition, recalling that
$\delta_\perp=\omp-\Omega\,$.
\begin{proposition}~\\
There exists a constant $C >0$ such that, for all $\ep\in ]0,1]$, for
all $\omp$, $\Omega$ s.t. $\dep\geq 1$ and for all  $ g \geq 0\,$, \beq
\label{precmean}
\inf_{||\Psi||=1} Q_{\Omega}^{per,N}(\Psi)= \omp +
m_A^N(\ep,\widehat g) \left( 1 -C r_A(\ep,\widehat g)  \right)\,, \eeq
with 
\beq\label{precmeana} 0 \leq r_A(\ep,\widehat g) \leq
g^{1/4}\delta_\perp^{-\frac
  18} \left(\frac{\dep +\omp}{\dep}\right)^\frac 14 \,
m_A^N(\ep,\widehat g )^{\frac 14} +    \,m_A(\ep,\widehat
g)\dep^{-1} \,. \eeq
\end{proposition}

\paragraph{Proof of the proposition}~\\ For simplicity, we make the
proof for  $\Omega=0$. Note also that
$$
1- Cr_A(\epsilon, \widehat g) \geq 0\,$$ by the lower bound. So we
 have only to prove \eqref{precmeana} under the additional condition
 that the right hand side of \eqref{precmeana}
 is less than some fixed $\alpha_0\,$. In any case, the estimate is only
 interesting in this case ! \\
   The proof does not depend on $N$ and  for
$\Omega$ not zero, we will make a remark at the end on how to adapt
it, using the diamagnetic inequality.

 The proof is inspired
by \cite{AB} where a reduction is made from a 3D to a 2D setting for
a fast rotation. We project a minimizer $\Psi$ onto
 $\psi_\perp\otimes L^2 (\mathbb R/N T \mathbb Z)$, and call
 $\psi_\perp(x,y)\,\xi(z)$ its projection:
 \beq\label{pro}
 \Psi(x,y,z)=\psi_\perp(x,y)\xi(z)+w(x,y,z)\hbox{ with }\int_{\rz^2}\psi_\perp(x,y)w(x,y,z)\, dxdy=0\,.\eeq
 The orthogonality condition implies in particular
\beq\label{gap}
 1=\int_{-\frac{NT}{2}}^{\frac{NT}{2}}|\xi(z)|^2\,dz\;+
\;\int_{\rz^2\times]- \frac{NT}{2},\frac {NT}{2}  [}|w(x,y,z)|^2\,dxdydz
\eeq
and we have the lower bound
\beq\label{lwgap}
\int_{-\frac{NT}{2}}^{\frac{NT}{2}} \Epsilon'_B (w(\cdot,\cdot,z))\;dz \geq
 (\dep +\omp)\, \int_{\rz^2\times ]- \frac{NT}{2},\frac {NT}{2}  [}|w(x,y,z)|^2\;dxdydz\,,
\eeq
with
$$
\Epsilon'_B(\psi) = \int_{\mathbb R^2}\left(
 \frac 12 |\nabla_{x,y}\psi (x,y) |^2 +
 \frac{\omp^2}{2}(x^2+y^2)\,|\psi(x,y)|^2\right)\, dxdy\,.
$$
 We compute the energy of $\Psi$ and use the orthogonality
  condition and the equation satisfied by $\psi_\perp$
 to find that all the cross terms
  disappear so that
  \begin{multline}\label{enint}
  Q^{N,per} (\Psi)= \omp \int_{-\frac {NT} 2}^{\frac {NT} 2} |\xi(z)|^2\, dz\,+\, \Epsilon^{N\,'}_A(\xi)\\
+\int_{\rz^2}\Epsilon^{N\,'}_A(w(x,y,\cdot))\ dxdy
  +\int_{-\frac {NT} 2}^{\frac {NT} 2} \Epsilon'_B (w(\cdot,\cdot,z))\ dz\\
  +g \int_{\mathbb R^2 \times ]- \frac{NT}{2},\frac {NT}{2}  [}|\Psi(x,y,z)|^4\,dxdydz\,,
\end{multline}
  where
$$
\Epsilon^{N\,'}_A(\phi)=\int_{-\frac {NT}2}^{\frac {NT} 2} \left(\frac12  |
 \phi'(z)|^2+W_\ep(z )|\phi|^2\right) \, dz\,.
$$From (\ref{gap}), (\ref{lwgap}) and  (\ref{enint}), we find \beq\label{lwbd}
Q^{N,per} (\Psi)\geq\omp+ \frac{\dep}{\dep +\omp}
  \int_{-\frac {NT}2}^{\frac{NT}2} \Epsilon'_B (w(\cdot,\cdot,z))\ dz +\int_{\rz^2}\Epsilon^{N\,'}_A(w(x,y,\cdot))\ dxdy\,.
\eeq
We use \eqref{lwbd}
 together with the upper bound \eqref{easyubA} and (\ref{lwgap}) to
 derive
  that
\beq \label{majw}
\int_{\rz^2\times]-\frac {NT}2, \frac{NT}2 [}|w(x,y,z)|^2\,dxdydz\leq
\frac{m_A^N(\ep,\widehat
 g)}{ \dep} \,.
\eeq
Note that the righthand side in \eqref{majw} is very small according
to Conditions \eqref{newAa} and \eqref{newAb}.\\
Note that \eqref{majw}  implies
\begin{equation}\label{minxi}
\int_{-\frac {NT}2}^{\frac {NT}2} |\xi(z)|^2 dz \geq  1-\frac{m_A^N(\ep,\widehat
 g)}{ \dep} \,.
\end{equation}
 Then, we get  also,
 \beq
\begin{array}{l}
\int_{\rz^2\times]-\frac {NT}2,\frac {NT}2[}|\nabla_{x,y} w(x,y,z)|^2 \,dxdydz\leq 2
\frac{\dep +\omp}{\dep} \; \frac{m_A^N(\ep,\widehat g)}{\omp},\\
\int_{\rz^2\times]-\frac {NT}2,\frac {NT}2[}|\partial_{z} w(x,y,z)|^2\,dxdydz \leq   2 \; m_A^N(\ep,\widehat g)\,.
\end{array}
\eeq
The proof of the  Sobolev embedding of $H^1(\mathbb R^3)$ in
 $L^6(\mathbb R^3)$
 gives (see for example \cite{Bre}, p.~164, line -1) for a general function $v$
 in  $H^1(\mathbb R^3)$
\beq
\|v\|_6\leq 4 \|\partial_x v\|_2^{1/3}\|\partial_y
v\|_2^{1/3}\|\partial_z v\|_2^{1/3}.
\eeq
Here $\|\;\cdot\|_p$ denotes the norm in $L^p(\mathbb R^3)$.\\
In our case, we are working
in $H^1(\mathbb R^2_{x,y} \times (\mathbb R_z/NT\mathbb Z))$. A
partition of unity in the $z$ variable allows us to extend this
estimate also this case, and we get, for another universal constant
$C$,
\begin{equation}
\|w\|_{6}
 \leq C_N  \|\partial_x w\|_2^{1/3}\|\partial_y
w\|_2^{1/3}\left(\|\partial_z w\|_2^2 + ||w||_2^2\right)^{1/6}\,,
\end{equation}
where this time $||\;\cdot\;||_p$ denotes the norm in
$L^p( \mathbb R^2_{x,y}\times ]-\frac {NT}2,\frac {NT} 2[)$.\\
So we obtain~:
\begin{equation}
\|w\|_{6}\leq {\tilde C}
m_A^N(\ep,\widehat g )^\frac 12 \left(\frac{\dep
  +\omp}{\dep}\right)^\frac 13\,.
\end{equation}
($C$, $\tilde C$ are $N$-dependent constants possibly changing from line
to line.)\\
Since by H\"older's Inequality,
$$
\|w\|_4\leq  \|w\|_2^{1/4}\|w\|_6^{3/4}\,,
$$
 we deduce that
\beq \label{wl4}
\|w\|_4\leq C \;m_A(\ep,\widehat g)^{\frac 12} \dep^{-\frac 18} \left(\frac{\dep
  +\omp}{\dep}\right)^\frac{1}{4} \,.
\eeq
We expand
$$
 |\Psi|^4=|\psi_\perp|^4|\xi|^4+2|\psi_\perp|^2 |\xi|^2 |w|^2 + 4(\Re(\psi_\perp \xi
\overline w)+\frac12 |w|^2)^2
  +4 |\psi_\perp|^2 |\xi|^2\Re(\psi_\perp  \xi \overline w)\,.
$$
  Since (\ref{enint}) implies that
  $$
\Epsilon^N(\Psi)\geq \omp+\Epsilon^{N}_A(\xi)- 4 g
 \int_{\mathbb R^2\times ]-\frac {NT}2,\frac {NT}2[} |\psi_\perp(x,y)|^3|\xi(z)|^3|w(x,y,z)|\,dxdydz \,,$$
 in order to get the lower bound, we just need to prove that the last
 term  is a perturbation to $\Epsilon^{N}_A(\xi)$. \\
We can do the following estimates
 $$
\begin{array}{l}
g   \int |\psi_\perp(x,y)|^3|\xi(z)|^3|w(x,y,z)|\,dxdydz\quad \\
\quad\quad\quad \leq c_0 g\omp^\frac 34  (\int
|\psi_\perp(x,y)|^4\,dxdy)^{\frac 34} \, (\int |\xi(z)|^4
 dz)^{\frac 34} \,\|w\|_4\\
\quad\quad\quad  \leq c_1
 g^{1/4}(\Epsilon^N_A(\xi))^{3/4}\|w\|_4\\
\quad\quad \quad \leq  c_2 g^{1/4}\dep^{-\frac 18}  \left(\frac{\dep
  +\omp}{\dep}\right)^\frac{1}{4} \,
m_A^N(\ep,\widehat g)^{\frac 12}
 (\Epsilon_A^N(\xi))^{3/4}\\
\quad\quad \quad \leq  c_3 g^{1/4}\dep^{-\frac
18} \left(\frac{\dep
  +\omp}{\dep}\right)^\frac{1}{4}\,m_A^N(\ep,\widehat g)^{\frac 14}
 \left(1 + C \,m_A^N(\ep,\widehat g)\dep^{-1}\right)
 \; \Epsilon_A^N(\xi)\,.
\end{array}
$$
Here to get the last line, we have used  the lower bound
 $$
\Epsilon_A^N(\xi) \geq m_A^N(\ep,\widehat g)\, ||\xi||_2^4 \,,
$$
 and \eqref{minxi}.\\
 This leads to
$$
\Epsilon^N(\Psi)\geq \omp+\Epsilon_A^N(\xi) \left(1 - C\, g^{1/4}\dep^{-\frac
  18}  \left(\frac{\dep
  +\omp}{\dep}\right)^\frac{1}{4} \,m_A^N(\ep,\widehat g )^{\frac 14}
-
 C \,m_A^N(\ep,\widehat g) \dep^{-1} \right) ,
$$
and then to \eqref{precmean}.\\

 \begin{remark}~\\
 In the case with rotation $\Omega$, the proof is the same if we replace $\Epsilon'_B$ by $\Epsilon'_{B,\Omega}$ defined by \beq
\Epsilon'_{B,\Omega} (\psi)
 = \int_{\mathbb R^2}
\left( \frac 12 |\nabla_{x,y} \psi - i \Omega r^\perp \psi|^2
 + \frac 12 (\omp^2 - \Omega^2) r^2|\psi|^2 \right)
 \, dx dy\,.
\eeq We also use the diamagnetic inequality \beq\label{diamag} \int
|\nabla |w| (x,y)|^2\, dx dy \leq \int |\left(\nabla w - i \Omega
r^\perp w\right)(x,y)|^2 \, dx dy \eeq which provides the Sobolev
injections.
\end{remark}

\begin{remark}~\\
 Here, we have not proved that the minimizer of $\Epsilon$ behaves almost
like the ground state in $x,y$ times a function of $\xi$ which
minimizes $\Epsilon_A$.  We are just able (see \eqref{majw}) to
prove that the minimizer is close to its projection (in some $L^2$
or $L^4$ norm). When $N=1$, this can  be improved under the stronger condition
\eqref{condbizb}.  We first observe (note that \eqref{lwbd} is still
true with
 the addition of $\Epsilon'_A(\xi)$ on the right hand side) that
\beq \Epsilon'_A(\xi) \leq m_A(\ep,\widehat g)\,. \eeq Using
\eqref{minxi}, assuming $\frac{m_A}{\dep} <1$, this leads  to \beq
\Epsilon'_A(\xi) \leq m_A (\ep,\widehat g) (1-
\frac{m_A(\ep,\widehat g)}{\dep})^{-1} ||\xi||^2 \eeq
We will show in Subsection \ref{sssB12} (see \eqref{compa}) how to
proceed in order
 to show that $\xi$ is close to the ground state
 $\phi_1(z)$ of $H_z^{per}$.\\
This can allow to improve the information given in Theorem
\ref{ThQL}.
\end{remark}



\section{The $1D$ periodic model~: estimates for $m_A^N$}\label{s1D}

The aim of this section is to analyze
$m_A^N$. We note that rough estimates were already given for the weak
interaction case which were enough for the justification of the model
 but  the corresponding rough estimates needed for the Thomas-Fermi justification
 will be obtained in this section. We will then look at
 accurate estimates for $m_A^N$, which will be established under
 stronger hypotheses. We will end the section by the discussion of the case
 $N>1$, which finally leads to the introduction of the DNLS model
 for the Weak Interaction case.

\subsection{ Universal estimates}\label{sssB11}
We consider the  one dimensional situation and a $T$- periodic
potential
 $W$, which could  be for example
$W(z)=  (\sin \pi  z) ^2/\ep^2$. We consider the problem of
minimizing on $L^2(\mathbb R/T \mathbb R)$
 the functional
\begin{equation}\label{B.1}
\psi \mapsto \mathcal G (\psi) = \frac 12\int_{-\frac T2}^{\frac T2}
|\psi'(z)|^2 \,dz + \int_{-\frac T2}^{\frac T2} W (z)|\psi(z)|^2\,
dz + \widehat g \int_{-\frac T2}^{\frac T2} |\psi(z)|^4 \, dz \,,
\end{equation}
over $||\psi||_{L^2}=1$.\\
We are interested in the control of the minimum of the functional
 and will simply prove
\begin{lemma}\label{approx1}~\\
If $\widehat g\geq 0$, then
\begin{equation}
m(\widehat g):=\inf_{||\psi||_{L^2}=1} \mathcal G (\psi)
 = \lambda_1 + \widehat g \int_{-\frac T2}^{+\frac T2}  |\phi_1(z)|^4 \,dz +
 o(\widehat g
 )\,,
\end{equation}
where $(\lambda_1, \phi_1)$ is the spectral pair
 of $-\frac 12\,\frac{d^2}{dz^2} + W(z)$ corresponding to the ground state
 energy (with $||\phi_1||^2=1$).
\end{lemma}

\begin{proof}~\\
It is clear that
\begin{equation}
\lambda_1 \leq  m(\widehat g) \leq  \lambda_1 + \widehat g
\int_{-\frac T2}^{\frac T2} |\phi_1(z)|^4 dz\,,
\end{equation}
so the question is now to improve the lower bound.\\
One could of course think of applying bifurcation theory
 but this gives only a local result and we need  in any case a global
 estimate
 for showing that the global minimizer of $\mathcal G$
 is closed to $\phi_1$ as $\widehat g$ is small.\\
Let $\phi_{min}$ be a minimizer of $\mathcal G$, then we know that
\begin{equation}
\frac 12  \int_{-\frac T2}^{\frac T2} |\phi_{min}'|^2 dz +
\int_{-\frac T2}^{\frac T2} W(z) |\phi_{min}(z)|^2\, dz \leq
\lambda_1+ \widehat g \int_{-\frac T2}^{\frac T2} |\phi_1(z)|^4\,
dz\,.
\end{equation}
So $\phi_{min}$ plays the role of a quasimode (or approximate
 eigenfunction) for
 $-\frac 12 \frac{d^2}{dz^2} + W(z)$.\\
A rather standard theorem in perturbation theory (we can write
 $\phi_{min} =\alpha \phi_1 + u^{\perp}$), gives first
$$
1-|\alpha|^2=|| u^{\perp}||^2 \leq \widehat g\frac{ \int_{-\frac
T2}^{\frac T2} |\phi_1(z)|^4 \, dz}{\lambda_2 -
  \lambda_1}\,,
$$
 and then the existence of a complex number $c$ of modulus $1$  such that
\begin{equation}\label{contapp}
||\phi_{min} - c \phi_1||^2_{L^2} \leq 2 \widehat g\frac{
\int_{-\frac T2}^{\frac T2} |\phi_1(z)|^4\,  dz}{\lambda_2 -
\lambda_1}\,.
\end{equation}
Here $\lambda_2$ denotes the second eigenvalue of or Hamiltonian.\\
Of course, the estimate is only interesting if
\begin{equation}\label{condrest}
2 \widehat g \frac{ \int_{-\frac T2}^{\frac T2} |\phi_1(z)|^4\,
dz}{\lambda_2 - \lambda_1}<1\,.
\end{equation}

We can now write
\begin{equation}
\begin{array}{ll}
m(\widehat g) & \geq \lambda_1 + \widehat g \int_{-\frac T2}^{\frac T2} |\phi_{min}(z)|^4\, dz\\
& \geq \lambda_1 + \widehat g \int_{-\frac T2}^{\frac T2}
|\phi_1(z)|^4\, dz
 - 4  \widehat g ||\phi_{min}-c \phi_1||_2 ||\phi_1||_{6}^3 \,.
\end{array}
\end{equation}
For the last estimate, we develop $|\phi_{min}|^4$ in the following
way
\begin{equation}
\begin{array}{ll}
|\phi_{min}|^4 &= |c\phi_1 +\phi_{min} - c\phi_1|^4 \\&
 \geq |\phi_1|^4 + 2 |\phi_1|^2 |\phi_{min} - c \phi_1|^2 - 4
 |\phi_1|^3 |\phi_{min} - c \phi_1|\,.\end{array}
\end{equation}
~From this inequality we get
\begin{equation}\label{min1}
|\phi_{min}|^4  \geq |\phi_1|^4  - 4 |\phi_1|^3 |\phi_{min} - c
\phi_1|\,.
\end{equation}
It just remains to control $||\phi_1||_6$ uniformly
 with respect to $\widehat g$, which can be deduced of the uniform control
 of the norm of $\phi_1$ in $L^6$.
\end{proof}

One can actually be more precise on what we have claimed in Lemma~\ref{approx1}.
\begin{lemma}\label{approx2}~\\
If $\widehat g\geq 0$, then
\beq \label{remain1}
 m(\widehat g) \geq \lambda_1
 + \widehat g ||\phi_1||_4^4  - 2^{\frac 52} \widehat g^\frac 32 ||\phi_1||_6^3
||\phi_1||_4^2 (\lambda_2-\lambda_1)^{-\frac 12}\,. \eeq This
estimate is interesting if \beq
 \widehat g
< \frac {1}{32}  ||\phi_1||_4^4 ||\phi_1||_6^{-6}
(\lambda_2-\lambda_1)\,. \eeq
\end{lemma}
\begin{remark}~\\
Everything being universal, one can of course replace $T$ by $NT$
 in the description.
\end{remark}

\subsection{Semi-classical results in the Weak Interaction case~:
 $N=1$}\label{sssB12}
We first recall that using \eqref{normel4a} we have,
 under Condition \eqref{AsscasAb},
 the rough control
\begin{equation}\label{proofapproxmAN=1}
\frac{1}{C \epsilon} \leq \lambda_{1,z} \leq m_A(\ep,\widehat g)\leq  \lambda_{1,z}
 + \widehat g \int_{-\frac T2}^{\frac T2}
|\phi_1(z)|^4\, dz \leq \frac{C}{\ep}\,,
\end{equation}
which leads to \eqref{approxma} for $N=1$ and was sufficient for the
justification of the longitudinal  model A.\\

Let us now show that under stronger assumptions
 one can have a more accurate asymptotics including the main
 contribution of the non-linear interaction.

\begin{proposition}~\\
Under Assumption  \eqref{condbizb}, $m_A$ admits the following asymptotics~:
\begin{equation}\label{estimategs} m_A(\ep, \widehat g)
=\lambda_1^{har}(\ep)  + \pi^{-\frac 1 2}{\bf w}''(0)^{\frac 14} \widehat
g \ep^{-\frac12} + c_0+ \mathcal O(\ep)+ \mathcal O (\widehat
g^\frac32 \ep^{-\frac 14})\,.
\end{equation}
\end{proposition}

\begin{proof}~\\
Indeed, $\lambda_1$ and $\lambda_1-\lambda_2$ are of order $\frac 1 \ep$,
 and by
\eqref{normel4a} and \eqref{contapp},  we get
\beq \label{compa} ||\phi_{min} - c \phi_1||^2_{L^2}\leq C \widehat
g \ep^\frac 12\,. \eeq
 Using the harmonic approximation, the term  $||\phi_1||_{6}$
  is of order $\ep^{-\frac 16}$
 and the remainder appearing in \eqref{remain1}
 is of order $\widehat g^\frac32 \ep^{-\frac 14}$.
 Altogether we get for the energy
\begin{equation}\label{mAimp}
m_A(\ep, \widehat g) =  \lambda_{1,z} + \widehat g \int_{-\frac
T2}^{\frac T2} |\phi_1(z)|^4 \,dz\;+  \mathcal O (\widehat g^\frac32
\ep^{-\frac 14})\,.
\end{equation}
Using \eqref{normel4a}, we obtain \eqref{estimategs}. This
asymptotics  becomes interesting in the semi-classical regime if
 \eqref{condbizb} holds.
\end{proof}
\begin{remark}~\\
Exponentially small effects will be discussed in Section
\ref{stunneling}.
\end{remark}


\subsection{Semi-classical analysis  in a Thomas-Fermi regime~: case $
  N=1$.}\label{ssTFr}
\subsubsection{Main results}\label{ssstfmr}
In this subsection,   we first give  the rough estimate leading to
 \eqref{approxmATF} for $N=1$. Recall that
 $\widehat g= \frac 1 \pi\,g\omp$, but $\widehat g$ and $\ep$  are taken as
 independent parameters.

\begin{proposition}\label{proproughtf}~\\
If for some $c >0$,
\beq\label{assthbww}
\widehat g \ep^2 \leq c\,,
\eeq
and if
\beq \label{assthf}
\widehat g \ep^\frac 12 >>1\,,
\eeq
then there exist $C$ and $\ep_0$ such that
\begin{equation}\label{rhf1}
\frac 1C\,\widehat g^\frac 23 \ep^{-\frac 23} \leq m_A(\ep,\widehat g) \leq C \widehat g^\frac 23 \ep^{-\frac 23}\,,\; \forall
\ep\in ]0,\ep_0]\,.
\end{equation}
\end{proposition}

We will also get the following accurate estimate~:

\begin{proposition}\label{proptfreg}~\\
 If
\beq\label{assthfb} \widehat g \ep^2  <<1\,, \eeq
and \eqref{assthf}
are  satisfied , then
\begin{equation}\label{assmA}
m_A (\ep,\widehat g) = 2^{-\frac 43}3^{\frac 53}5^{-1} {\bf w}''(0)^{\frac
23} \widehat g^{\frac 23} \ep^{-\frac 23}
 \left(1+ \Og (\widehat g^{-\frac 23} \ep ^{-\frac 13})\right) \,.
\end{equation}
\end{proposition}
 The new asusmption is \eqref{assthfb}, which is stronger than \eqref{assthbww}.

\subsubsection{The harmonic functional on $\mathbb R$}\label{ssshfR}
 Let us start with
 the case of a harmonic potential
 $W_\ep(z)= \gamma \frac{z^2}{2 \ep^2}$ on $\mathbb R$, with $\gamma
 >0$,  and consider
 the problem of minimizing
\beq
q^{Hr,T} (u) = \frac 12  \int_{-\frac T2}^{\frac T2} u'(t)^2\, dt
+ \frac{\gamma}{2 \ep^2}  \int_{-\frac T2}^{\frac T2} t^2u(t)^2\,dt
 + \widehat g  \int_{-\frac T2}^{\frac T2} u(t)^4\,dt
\eeq
over the $u$'s in the form domain of $q^{Hr,T}$ such that $||u||^2=1$.\\
We denote by $m_A^{Hr,T}$ the infimum of the functional.
Actually there are two approximating  `` harmonic ''
functionals  of interest corresponding to $T$ finite and to $T=+\infty$.
An interesting  point is that,
 for $T$ large enough,
 the minimizers of these
two functionals are the same as we will see below. But let us start
with
 the case $T=+\infty$.
\begin{lemma}~\\
 If \eqref{assthf}
holds, then
\begin{equation}\label{5.22}
m_A^{Hr,+\infty} (\ep,\widehat g) = 2^{-\frac 43}3^{\frac 53}5^{-1}
\gamma ^{\frac 23} \widehat g^{\frac 23} \ep^{-\frac 23} \left(1+ \Og
(\widehat g^{-\frac 23} \ep ^{-\frac 13})\right) \,.
\end{equation}
\end{lemma}
The proof is very standard (see for example \cite{BBH}, \cite{Af} or  \cite{CorR-DY1}
 which treat the $(2D)$-case).
The analysis is done through a dilation. We look for an  $L^2$-normalized
 test function $\phi$ in the form
\beq \label{dilat}
\phi(z) = \rho^\frac 12 v(\rho z)\,,
\eeq
with $\rho$ and $v$ to be determined.\\
 The $1-D$ energy of $\phi$ becomes
\beq \label{renorm}
\frac 12 \rho^2 \int_{\mathbb R} v'(t)^2 dt + \rho^{-2} \ep_\gamma^{-2}
\int_{\mathbb R} t^2 v(t)^2 dt
 + \widehat g \rho \int_{\mathbb R} v(t)^4 dt\,,
\eeq
with
$$
\ep_\gamma = \ep/\sqrt{\frac 12 \gamma}\,.
$$
This leads to  choose $\rho=\rho_\gamma $ such that
$$
 \rho^{-3} = \widehat g
\ep_\gamma^2\,,$$
 hence
 \beq \label{defrho}
\rho_\gamma  = \ep_\gamma^{-\frac 23} \widehat g^{-\frac 13}\,,
\eeq
 and the energy of this model becomes
\beq
\widehat g^{\frac 23} \ep^{-\frac 23} \left( q_{TF}(v) + \frac
12(\ep_\gamma^\frac 12 \hat g)^{-\frac 43} \int_{\mathbb R}
v'(t)^2\,dt
\right)
\eeq
with
\begin{equation}
q_{TF}(v):= \int_{\mathbb R}  t^2 v(t)^2\,dt +
\int_{\mathbb R} v(t)^4\, dt\,.
\end{equation}
This is
asymptotically of the order of $
\widehat g^{\frac 23} \ep^{-\frac 23}\;$ and
Condition \eqref{assthf} is just the condition  that the kinetic term is negligeable
 in the computation of the energy.\\
Let us recall the details of this asymptotics  for completeness. We have
first to minimize over $L^2$-normalized
 $v$ the approximating functional  $q_{TF}$.
The minimizer $ v_{min} (t)$ of $q_{TF}$ is determined by the
equations
\begin{equation}
t^2 v(t) + 2 v(t)^3 = \lambda v(t)\,,\; \int_{\mathbb R} v(t)^2
dt=1\,.
\end{equation}
We get
\begin{equation}
v_{min}(t) = 2^{-\frac 12} (\lambda -t^2)_+^\frac 12 \,,
\eeq
with
\beq
\lambda = \left( \frac 32\right)^\frac 23\,,
\eeq
and for $x\in \mathbb R$, 
$$
(x)_+=\max(x,0)\,.
$$
The corresponding TF-energy is \beq e_{TF}:=\int_{\mathbb R} (t^2
v_{min}(t)^2 + v_{min}(t)^4)\,dt =\frac 25  \lambda^\frac 52 \,.\eeq
Unfortunately, this minimizer is not in $H^1(\mathbb R)$
 and can not be used directly for our initial rescaled functional
$$
q_{TF}^{\sigma}(v) = q_{TF}(v) + \sigma \int_{\mathbb R}
v'(t)^2\,dt\,,
$$
with
$$
\sigma = \widehat g^{-\frac 43} \ep_\gamma^{-\frac 23}\,.
$$
Here we recall that \eqref{assthf} implies
 $$0\leq \sigma <<1\,.
$$
 So we need
 to regularize this minimizer to have an upperbound
 for the energy of our ``harmonic'' functional which is good as
 $\sigma \ar 0$.\\

This can be done in the following way (see for example \cite{Af} and
references therein).

We introduce \beq \gamma(s) =\left\{
\begin{array}{ll}
\sqrt{s} \,,\; \mbox{ if } s > \sigma^\frac 13\\
s \sigma^{-\frac 16} \,,\; \mbox{ if } s < \sigma^\frac 13
\end{array}\right.
\eeq Let us consider the  function
$$
\hat v_\sigma (t) = \gamma (v_{min}(t)^2)\,.
$$
We get that $\hat v_\sigma$ belongs to $H^1$ and satisfies
$$
\int |\hat v_\sigma(t)|^2 \,dt =1 - r(\sigma)\,,
$$
with
$$ r(\sigma) = \Og (\sigma^\frac 23)\,.
$$

More precisely the (positive) remainder $r(\sigma)$  is
$$
r(\sigma)= \int_{v_{min}(t)^2 < \sigma^\frac 13} (-|v_{min}(t)|^2 +
 \sigma^{-\frac 13} |v_{min}(t)|^4)\, dt \,.
$$
Let us now consider as a test function
$$
v_\sigma := \hat v_\sigma /\|\hat v_\sigma\|\,.
$$
Then we have
$$
v_\sigma = (1 +\frac 12 r(\sigma) +\Og( \sigma^\frac 43)) \hat
v_\sigma\,.
$$
So  we get
$$
q_{TF}^\sigma (v_\sigma)= q_{TF}(v_{min}) + \Og (\sigma \ln \frac
1\sigma)\; .
$$

So we have the upper-bound in the statement \eqref{5.22} (actually with a
better remainder term)  of the lemma. The lower
bound in \eqref{5.22}  is immediate
 because the kinetic term is positive.
\subsubsection{The harmonic functional on $]-\frac T2,\frac T2[$}\label{ssshfT}
We consider now 
the case of the interval 
 and have the following Lemma~:
\begin{lemma}\label{HarmT}~\\
Under Assumption \eqref{assthf}, there exists $C>0$ such that
\beq
m_A^{har, T}(\ep,\widehat g)\geq \frac 1C\, \widehat g^{\frac 23} \ep^{-\frac 23}\,.
\eeq
\end{lemma}
The proof is based on  the same method as in the previous
subsubsection.  It is easy to see that the minimizers coincide if
\beq
\frac {\rho_\gamma  T }{2} > \lambda^\frac 12 \,,
\eeq
that is
\begin{equation}\label{condsurT}
T > \widehat g^{\frac 13} \ep_\gamma^\frac 23 \left(\frac 32\right)^{\frac 13}\,.
\end{equation}

If \eqref{condsurT} is not satisfied, we can still have a lower bound
for the infimum of the functional. The renormalized functional reads
\beq \label{renorma}
q^{ren,T}(v):=\rho^2 \int_{\frac {\rho T }{2}} ^{\frac {\rho T }{2}} v'(t)^2 dt + \rho^{-2} \ep_\gamma^{-2}
\int_{\frac {\rho T }{2}}^{\frac {\rho T }{2}}  t^2 v(t)^2 dt
 + \widehat g \rho \int_{\frac {\rho T }{2}}^{\frac {\rho T }{2}}  v(t)^4 dt\,,
\eeq
which satisfies
$$
q^{ren,T} (v)\geq \widehat g \rho \left(\int_{\frac {\rho T }{2}}^{\frac
  {\rho T }{2}}  v(t)^4 dt\right)\,.
$$
Using the H\"older inequality, we obtain, if $|| v||_2=1$,
$$
q^{ren,T}(v) \geq (\widehat g \rho ) (\rho T)^{-1}\,,
$$
and using our assumption, we obtain
\begin{equation}
q^{ren,T}(v) \geq \frac 12 \lambda^{-\frac 12}
 ( \widehat g \rho ) \geq  \frac 1C \widehat g^\frac 23 \ep^{-\frac 23}\,,
\eeq
if $||v||_2=1$.

We then immediatly obtain Lemma \ref{HarmT}.

\subsubsection{Relevance of the ``harmonic functional'' for rough
  bounds}\label{sssRhfrough}
First we prove Proposition \ref{proproughtf}. We can proceed by direct comparison. Observing that we can find
$\alpha >0$ such that
$$
{\bf w}(z) \leq \alpha z^2\,,\; \forall z\in [-\frac T2,+\frac T2]\,,
$$
and
$$
\rho_\alpha T > 2 \lambda^\frac 12\,.
$$
Here, we use \eqref{assthbww} and
$$
\rho_\alpha = c_0\alpha^{\frac 13} (\ep^{-\frac 23} \widehat g^{-\frac
  13}) \geq c_0\alpha^{\frac 13} c^{-\frac 13}\,.
$$
We can then use the asymptotic estimate \eqref{5.22} with
$\gamma=\alpha$
 to  get the upper bound in \eqref{rhf1}.\\

Using now Assumption \eqref{assw}, we can also find $\hat \alpha$
 such that
$$
{\bf w}(z) \geq \hat \alpha z^2\,,\; \forall z\in [-\frac T2,+\frac T2]\,,
$$
This leads, using our analysis of $q^{TF}$ in the harmonic case
 to the lower bound in \eqref{rhf1}.

\subsubsection{Relevance of the ``harmonic functional' for the asymptotic
 behavior}\label{ssshfacc}

In order to have a  better localized minimizer, we should assume
 that $\rho \rightarrow +\infty$ and this corresponds to replacing
 Assumption  \eqref{assthbww} by the stronger Assumption
 \eqref{assthfb}.

Moreover, we have to verify that under this assumption
 the ``harmonic approximation'' is valid
 for this energy computation. For this, we should analyze the
 localization of the minimizer.
Assuming that such a localized minimizer exists (minimize the
functional
 $v\mapsto \int (z^2 v(z)^2 + v(z)^4)\, dz$),
 we can also get an upperbound of $m_A$ by using a harmonic
 approximation
 and a lower bound of the same order.

For the lower bound, we have just to analyze (forgetting the
positive
 kinetic term) the infimum
 of the functional
$$
\phi \mapsto \int_{-\frac T2}^{\frac T2} \left(\frac{{\bf w}(z)}{\ep^2}
\phi^2 + \widehat g \phi^4\right) \, dz\,.
$$
As in the other case, a minimizer (over the $L^2$-normalized
$\phi$'s),
 should satisfy, for some $\mu>0$,  the Euler-Lagrange equation
$$
\frac{{\bf w}(z)}{\ep^2} \phi(z) + 2 \widehat g \phi(z)^3 = \mu \phi(z)\,,
$$
where $\mu$ will be determined by the $L^2$ normalization over
$]-\frac T2, \frac T2[$.
\\
We find
\beq\label{mintf}
\phi(z)=\frac {1}{\sqrt{ 2\widehat g} }  \left(\mu -
\frac{{\bf w}(z)}{\ep^2}\right)^\frac 12 _+\,.
\eeq
with
\beq \label{tfnorm}
\frac {1}{2 \widehat g}\, \int  (\mu - \frac{{\bf w}(z)}{\ep^2}) _+ dz
=1\,.
\eeq
But we know from the upperbound that $\mu$ is less than two times
the energy which is asymptotically lower than $m_A^{har} (\ep
   \widehat g)$.
In particular, if $\mu \ep^2$ is small, it is easy to estimate $\mu$
 using the harmonic approximation of $w$ at its minimum.\\
It remains to verify the behavior of $\mu \ep^2$. We find
$$
\mu \ep^2 \leq C \widehat g^\frac 23 \ep^\frac 43\,.
$$
Not surprisingly, this shows that $\mu \ep^2$ is small as  $\rho \ar
+\infty$. So finally, we have obtained Proposition \ref{proptfreg}.

\subsection{The case $N>1$}~\\
We would like to extend our rough or accurate estimates
 for $m_A$ to the analogous estimates for $N>1$, keeping
 the same kind of assumptions.
\subsubsection{Universal control}
We now consider the functional over
$]-\frac{NT}{2},\frac{NT}{2}[$. Using the minimizer
 obtained for $N=1$ and extending it by periodicity, we get after
 renormalization, the general upper-bound
\begin{equation}\label{comparaison1N}
m_A^N (\ep,\widehat g)\leq m_A (\ep, \frac{\widehat g}{N})\,.
\end{equation}
~From this comparison, we obtain immediately the rough upper bounds
 in the WI case and in the TF case.\\

\subsubsection{Rough lower bounds}
In the WI case,  we always have, observing that $\lambda_{1,z}$
 is the ground state energy for any $N\in \mathbb N^*$,
\begin{equation}
\lambda_{1}^z \leq m_A^N (\ep,\widehat g )\,.
\end{equation}
Hence we obtain in full generality
\begin{proposition}\label{estigross}~\\
Under Condition (\ref{AsscasAa}), then, for any $N\geq 1$,
 we have
\begin{equation}
m_A^N (\ep,\widehat g) \approx \frac{1}{\epsilon}
\end{equation}
\end{proposition}

In the TF case, it remains to prove the lower bound which will be a consequence of
 the following inequality~:
\begin{equation}\label{tfN}
m_A^N (\ep,\widehat g) \geq \frac{1}{C N^2} \widehat g^{\frac 23}
\ep^{\frac 43}\,.
\end{equation}
We indeed observe that if $u_N$ is a normalized minimizer, then there
exists
 one interval $I_j:=]j \frac T2,(j+2) \frac T2[$ ($j\in \{-N,\dots,
    N-2\}$),  such that
$$
\int_{I_j} |u_N|^2\,dz\; \geq \frac 1N
$$
We can then write,  forgetting the kinetic term and translating
$I_j$ to $]-\frac T2,+\frac T2[$,
$$
\begin{array}{ll}
m_A^N (\ep,\widehat g) & \geq\ep^{-2} \int_{I_j} {\bf w}(z)\, |u_N|^2 \,dz
 + \widehat g \int_{I_j} |u_N|^4\, dz \\ &
\geq \inf (||u_N||^2,||u_N||^4) \inf_{||u||=1} \int_{-\frac
T2}^{+\frac
  T2} (W_\ep |u|^2 + \widehat g |u|^4)\,dz\,.
\end{array}
$$
Then we  can  combine
  the lower bound  obtained for $N=1$ and the inequality
 ${\bf w}(z)\geq \hat \alpha z^2$ to get \eqref{tfN}.
So we get finally that $m_A^N$ has the right order in the TF case.
\begin{proposition}\label{proofroughtfAN}~\\
Under Assumptions \eqref{assthbww} and  \eqref{assthf}, we have, for
any $N\geq 1$, \beq m_A^N (\ep,\widehat g) \approx \widehat g^{\frac
23} \ep^{\frac 43}\,. \eeq
\end{proposition}
This extends to general $N$  our former Proposition \ref{proproughtf}.
\subsubsection{Asymptotics}
We would like to give conditions under which the universal upperbound
 \eqref{comparaison1N} becomes actually asymptotically or exactly a
 lower bound.\\
\begin{proposition}\label{AsymptN}~\\
Under either Assumption \eqref{condbizb} or Assumptions \eqref{assthf} and \eqref{assthfb},
\begin{equation}\label{memeasympt}
m_A^N (\ep,\widehat g) \sim  m_A (\ep, \frac{\widehat g}{N})\,.
\end{equation}
\end{proposition}

\begin{proof}~\\
The upperbound was already obtained in \eqref{comparaison1N}. The proof
 of the lower bound is different in the two considered cases.\\

\paragraph{WI case.}
We will see later (in \eqref{asympassezprecise}) by a rough analysis
of the tunneling
 effect and the property that the infimum of the function
$$
\mathcal C^N\ni (c_0,c_2,\dots,c_{N-1}) \mapsto
 \sum_{j=0}^{N-1} |c_j|^4
$$
over $\sum_j |c_j|^2=1$ is attained when all the $|c_j|$'s are
equal~: \beq |c_j|=\frac {1}{\sqrt{N}}\,,\; \mbox{ for }
j=0,\dots,N-1\,, \eeq that, under Assumption \eqref{condbizb}, there
exist $C>0$, $\ep_0>0$ and $\alpha >0$
 such that
\beq
m_A^N (g,\ep) \geq  m_A (\frac{\widehat g}{N},\ep) - C  (\widehat g + 1)
\exp -\frac{\alpha}{\ep}\,,\;\forall \ep\in (0,\ep_0]\,.
\eeq

\paragraph{TF case.}
In this case we can for the lower bound forget the kinetic term
 and come back to the analysis of Subsubsection \ref{ssshfacc}, with $T$
 replaced by $NT$. Under Assumption \eqref{assthfb},
 we have seen in \eqref{mintf} that the minimizer $u_N$ is localized in the
 neighborhood of each minimum and $T$-periodic.\\
We can then write
$$
\begin{array}{ll}
\int_{-\frac{NT}{2}}^{\frac{NT}{2}}
 \left( \frac{{\bf w}}{{\ep}^2} |u_N|^2  + \widehat g |u_N|^4\right) \;dz
& = N \int_{-\frac{T}{2}}^{\frac{T}{2}}
 \left( \frac{{\bf w}}{\ep^2} |u_N|^2  + \widehat g |u_N|^4\right) \; dz\\
& = \int_{-\frac{T}{2}}^{\frac{T}{2}}
 \left( \frac{{\bf w}}{\ep^2} |\sqrt{N} u_N|^2  + \frac{\widehat g}{N}
|\sqrt{N} u_N|^4\right)\; dz\\&
\geq \inf_{||v||=1}  \int_{-\frac{T}{2}}^{\frac{T}{2}}
 \left( \frac{{\bf w}}{\ep^2} |v|^2  + \frac{\widehat g}{N}
|v|^4\right)\; dz\,.
\end{array}
$$
But under Assumptions \eqref{assthf} and \eqref{assthfb}, the last term in the
 inequality has same asymptotics as $m_A(\ep, \frac{\widehat g}{N})$
 and we are done.

\end{proof}
\begin{remark}~\\
This proposition leaves open the question of the equality in \eqref{memeasympt}.
\end{remark}


\section{Study of Case ($B$) : Justification of the transverse reduced model}\label{scasB}
\subsection{Main result}
 We have defined $\Epsilon_{B,\Omega}^N$ by
  \eqref{defEBN}-\eqref{defPsiN} and
 $m_{B,\Omega}^N$,
 the infimum of the energy by \eqref{defmB}.
In case B, the proof of the reduction does not depend  on whether
$N=1$ or $N>1$. The only difference is when looking at the rough or
accurate estimates of the reduced model. Note that only rough
estimates are used in the part concerning the justification of the
model.

The reduction is very similar to case A, and we will prove
\begin{theorem}~\\\label{theoRB}
If
\beq \label{RBa} (RBa)\quad \ep\, m_{B,\Omega}^N <<1\,,
\eeq
 and
\beq \label{RBb} (RBb)\quad g\, m_{B,\Omega}^N \,
\ep^\frac 12 <<1\,, \eeq
 then, as $\ep$ tends to $0$,
\beq
\inf_{||\Psi||=1} Q^{per,N}_{\Omega}(\Psi)  = \lambda_{1,z} + m_{B,\Omega}^N (1+ o(1))\,. \eeq
\end{theorem}
Then Theorems \ref{theoBWI} and \ref{theoBTF} follow from this
result and appropriate estimates on $m_{B,\Omega}^N$, as we will prove in
section \ref{s2D}, while the proof of Theorem~\ref{theoRB} is made
in Section \ref{srb}.

\subsection{Proof of Theorems \ref{theoBWI} and \ref{theoBTF}}\label{s2D}
The issue is to determine the magnitude of the infimum of the energy
of the transverse problem $m_{B,\Omega}^N$.
\subsubsection{Reduction to the case $N=1$}
As in Case $A$ it is immediate to see that
\beq\label{majcasB}
m_{B,\Omega}^N\leq
 m_{B,\Omega} (\frac{\widetilde g}{N},\omp)\,.
\eeq
If indeed $\psi_{min,N}$ was the $T$-periodic minimizer for \eqref{EB}
 with $\widetilde g_N =\frac{\widetilde g}{N}$, we get \eqref{majcasB}
 by  using
 \eqref {deunaN}, \eqref{deWaversphi} and taking
 $\psi_{j,\perp}= \frac{1}{\sqrt{N}} \psi_{min,N}$.\\

So it remains for the needed upper-bound to analyze the case $N=1$.
 This depends on
the magnitude of $\tilde g$ and  leads us  to consider two cases.
\subsubsection{The Weak Interaction regime : case $N=1$}
\begin{proposition}\label{majcasbWI}~\\
If (\ref{BWIa}) holds, then
\beq \label{majcasbWIa}
m_{B,\Omega}(\widetilde g,\omp)\leq C\omp\,.
\eeq
\end{proposition} Indeed,
(\ref{BWIa}) implies that $\tilde g$ is bounded and the test
function $\psi_\perp$ (which is independent of $\Omega$)
implies the proposition.

Therefore, if (\ref{BWIa}) and  (\ref{BWIb}) are
satisfied, then Theorem \ref{theoRB} holds and implies  Theorem~\ref{theoBWI}.

\subsubsection{The Thomas Fermi regime : case $N=1$}
We start with the case when $\Omega=0$. When  $\tilde g$ is not
 bounded, we can meet a Thomas-Fermi situation.
\begin{proposition}\label{estmB}~\\
If $\widetilde g \rightarrow +\infty$, the  function $m_B(\widetilde
g,\omp)$
 satisfies
\beq\label{ubb}m_B(\widetilde g,\omp)\sim c_{TF}
\omp\sqrt{\widetilde g
 }\,,
\eeq with \beq \label{calcc0} c_{TF}= \frac {\pi}{24} \lambda^3 =
3^{-1} 2^{\frac 32} \pi^{-\frac 12}\,. \eeq
\end{proposition}
Therefore, if (\ref{asscasBa}), (\ref{asscasBb}), (\ref{asscasBc})
are satisfied, then Theorem \ref{theoRB} implies Theorem
\ref{theoBTF}.

{\bf Proof.}~\\
 A rescaling in
$\sqrt{\sqrt{\widetilde g}/\omp}$  yields a new energy
$$
u\mapsto \frac \omp 2\; \int_{\mathbb R^2}\left( \frac
{1}{\sqrt{\widetilde g}}|\nabla u|^2+ \sqrt{\widetilde g} r^2 |u|^2
+ 2 \sqrt{\widetilde g} |u|^4\right) dx dy\,,
$$
which is of the type Thomas Fermi (that is kinetic energy can be
neglected) if \beq \label{widetildeggrand} \frac
{1}{\sqrt{\widetilde g}}<< \sqrt{\widetilde g}\,. \eeq This leads
then simply to
 the TF reduced functional
$$
u\mapsto (\omp \sqrt{\widetilde g}) \int_{\mathbb R^2}\left(\frac 12
r^2 |u|^2 + |u|^4\right)\; dx dy\,,
$$
whose infimum over the unit ball in $L^2(\mathbb R^2)$ is of order
$c_{TF}  (\omp \sqrt{\widetilde g})$, with $c_{TF} >0$ defined by~:
\begin{equation}\label{defc0tf}
c_{TF}= \inf_{||u||_2=1} \int_{\mathbb R^2} \left(\frac 12
r^2|u(x,y)|^2 +
  |u(x,y)|^4\right)\; dxdy\,.
\end{equation}
The minimizer exists and  is explicitly known as
$$
u_{min}(x,y) = \frac 12 (\lambda - r^2)_+^\frac 12\;\mbox{ with }
\lambda = 2^{\frac 32}\pi^{-\frac 12}\,.
$$
This leads to \eqref{calcc0}.

In addition, by  a computation similar to the one in  Subsection
\ref{ssTFr}, we obtain more precisely
\begin{lemma}~\\
There exists $c$ such that,
as $\widetilde g$ tends to  $+\infty$, \beq \frac{m_B}{\omp}  =
c_{TF} \, \sqrt{\widetilde g} +\frac{ c}{\sqrt{\widetilde g}} \ln
 \widetilde g + \mathcal O (\frac{ 1}{\sqrt{\widetilde g}})\,,
\eeq with $c_{TF}$ defined  in \eqref{defc0tf}.
\end{lemma}
\begin{remark}~\\
Note that we have the universal lower bound
\begin{equation}\label{newlbmB}
m_B(\widetilde g, \omp) \geq c_{TF} \, \omp \sqrt{\widetilde g}\,.
\end{equation}
This lower bound becomes better than the universal lower bound by $\omp$ as soon as \beq
\label{compoldnew} c_{TF}\,\sqrt{\widetilde g}  >1\,. \eeq
\end{remark}

\begin{remark}~\\
In the semi-classical regime,   conditions (BTFa) and (BTFc)
 in Theorem \ref{theoBTF} (take their product)
imply that this two-dimensional  energy is much smaller than
$1/\ep$, that is
\begin{equation}
\omp g^\frac 12 \ep^{-1/4}<<\ep^{-1}\,.
\end{equation}
\end{remark}

We now look at the case when $\Omega >0$. The previous proof, using
that the minimizer of the TF reduced functional in \eqref{defc0tf}
is radial, yields
\begin{proposition}~\\
There exists $C$ such that, as  $ \widetilde g \ar +\infty$,
\beq
m_{B,\Omega} (\widetilde g,\omp) \leq m_B(\widetilde g,\omp) + C \ln
\widetilde g \;\widetilde g^{-\frac 12}\,.
\eeq
\end{proposition}
This will be improved in \eqref{comprot}
 by a direct study of the minimizer of $\Epsilon_{B,\Omega}$.
\begin{remark}~\\
For a lower bound, we can use the TF reduced functional
$$
I_\Omega(u)= \omp\sqrt{\widetilde g} \;
 \int_{\mathbb R^2} \left(\frac 12 (1- \Omega^2/\omp^2) r^2 |u|^2
 + |u|^4\right)\; dx dy
$$ whose minimum is explicit~: $$
\inf _{||u||=1} I_\Omega(u) =\omp\sqrt{\tilde g}e_{TF}\sqrt
{\frac 12 (1- \Omega^2/\omp^2)}\,.
$$
 Thus we get that, if there exists
$\beta \in [0,1[ $ such that \beq\label{controleOmega}
 0\leq \Omega/\omp \leq \beta\,,
\eeq
then, as  $\widetilde g\ar +\infty$,
\beq\label{asympcasBOmegaN=1}
m_{B,\Omega}(\widetilde g,\omp)\approx \omp \sqrt{\widetilde g}\,.
\eeq
The uniformity of the approximation depends on $\beta$.

In fact, if one wants a more precise expansion of the energy,
one can use the ground state $\rho$ of
 $I_\Omega$ to split the energy $\Epsilon_{B,\Omega} (u)$. Indeed
 the Euler Lagrange equation for $\rho$ multiplied by $(1-|u|^2)$
  for any function $u$ yields the identity (see \cite{Af})
  $$ \Epsilon_{B,\Omega} (u)=I_{\Omega}(\rho)+\int\rho^2|\nabla
v-i\Omega\times r v|^2+\tilde g\rho^4
  (1-|v|^2)^2$$ where $v=u/\rho$. Thus, $I_\Omega$ always provides a lower
bound
   with an inverted parabola profile as soon as we are
    in a TF situation. The second part of the energy has
    the vortex contribution which is of lower order when
    $\Omega/\omp<<1$. More precisely, the first vortex is observed
     for a  velocity $\Omega$ of order $\omp \ln\tilde g/\sqrt{\tilde g}$.
     When $\Omega$ increases and becomes at most like $\beta\omp$ with
$\beta<1$,
      the two parts of the energy $I(\rho)$ and the rest become
       of similar magnitude. In the limit, $\Omega\to\omp$, there
       are a lot of vortices and the description can be made with
    the lowest Landau levels sets of states. The leading order term
     of the energy is the first eigenvalue of $-(\nabla -i\Omega\times r)^2$
     which is equal to $\Omega$.
\end{remark}

 \subsection{Proof of Theorem \ref{theoRB}}\label{srb}
We recall that
we have the universal upperbound \eqref{univubB}.
 The lower bound follows from the  following proposition and the
 fact that there exists $c>0$ such that $$
\delta_z^N \sim c/\epsilon\,,
$$
as $\epsilon$ tends to $0$.
\begin{proposition}\label{Prop6.9}~\\
 There exists a universal constant $C>0$ such that
\begin{equation}\label{minB}
\inf_{||\Psi||=1} Q_{\Omega}^{per,N}(\Psi)=  \lambda_{1,z} +
 m_{B,\Omega}^N \;( 1 -C  r_B^N) \,.
\end{equation}
with \beq
 0\leq r_B^N \leq   \,m_{B,\Omega}^N (\delta_z^N)^{-1}
+   g^\frac 14
(\delta_z^N)^{-\frac 18}( m_{B,\Omega}^N) ^\frac 14
 (1+\frac{\lambda_{1,z}}{\delta_z^N})^\frac 18  \,.
\eeq
\end{proposition}
Before giving the detailed proof, let us shortly sketch  the case
$N=1$. The proof is indeed essentially the same as for Case $(A)$.
One has just to exchange the role
 of $(A)$ and $(B)$.
$m_A$ should be replaced by $m_B$, $\omp$ by $\lambda_{1,z}$
 and $\delta_z$ by
$\delta_\perp$. Note that in the two models $(A)$ and $(B)$  the
ratio
 between the ground state energy and the splitting is of order $1$. We
 have indeed
$$
\delta_\perp =\omp\mbox{ and, in the semi-classical regime, }
\frac{\lambda_{1,z}}{\delta_z} \approx 1\,.
$$

\begin{proof}~\\
We start from a minimizer $\Psi$ and first write
\beq \label{app1}
\Psi = \Pi_N \Psi + w
\end{equation}
where $\Pi_N$ is the projection relative to the first $N$
eigenfunctions
 of $H_z$ introduced in \eqref{defPiN}. We have 
\beq\label{app2}
\Pi_N w=0\,,
\eeq
 and
\beq\label{app3}  ||w||^2 + ||\Pi_N \Psi||^2 =1\,. \eeq We have the
lower bound \beq \label{app4} \int_{\mathbb R^2_{x,y}}
\Epsilon'_A(w) \,dxdy\geq \lambda_{N+1,z} \int_{\mathbb R^2 \times
]-\frac{NT}{2},+\frac{NT}{2}[}|w(x,y,z)|^2 dxdydz\,, \eeq with \beq
\label{app5} \Epsilon'_A(\phi):=
\int_{-\frac{NT}{2}}^{\frac{NT}{2}}\left( \frac 12 \phi'(z)^2 +
 \frac {1}{\ep^2} {\bf w}(z) \phi(z)^2\right) \;dz\,.
\eeq
We now rewrite the energy in the form
\beq \label{app6}
Q^{per,N}_\Omega(\Psi) = \int_{-\frac{NT}{2}}^{\frac{NT}{2}}
\Epsilon'_{B,\Omega}(\Psi)\,dz
 + \int_{\mathbb R^2_{x,y}}\Epsilon'_A(\Pi_N\Psi) \,dxdy
 + \int_{\mathbb R^2_{x,y}}\Epsilon'_A(w) \,dxdy
 + I_N(\Psi)\,,
\eeq
with
\beq \label{app7}
I_N(\Psi)= g \int |\Psi|^4 dxdy dz\,,
\eeq
and
\beq \label{app8}
\Epsilon'_{B,\Omega} (\psi)= \int_{\mathbb R^2_{x,y}}
\left( \frac 12 |\nabla_{x,y} \psi -i
\Omega r_\perp\psi|^2 + \frac 12 (\omp^2-\Omega^2)
r^2|\psi|^2\right)\,dxdy\,,
\eeq
with
$r_\perp=(-y,x)$.

We note that $I_N\geq 0$ and that
\beq\label{app9}
\Epsilon'_{B,\Omega}(\psi) \geq \omp ||\psi||^2\,.
\eeq
We first get the control of $||w||^2$. Having in mind \eqref{univubB},
we obtain
\beq\label{app10}
\begin{array}{ll}
\lambda_{1,z} + m_{B,\Omega}^N  &\geq
 Q^{per,N}_\Omega(\Psi) \\ & \geq \omp  + \lambda_{N+1,z} ||w||^2 +
\lambda_{1,z} ||\Pi_N\Psi||^2
\end{array}
\eeq and this implies \beq \label{app11} ||w||^2 \leq
\frac{m_{B,\Omega}^N}{\delta_z^N}\,. \eeq The right hand side in
\eqref{app11} is small according to \eqref{RBa}. Note also that we
have immediately from \eqref{app3}, \beq \label{app12} ||\Pi_N
\Psi||^2 \geq 1 - \frac{m_{B,\Omega}^N}{\delta_z^N}\,. \eeq We now
have  to control the derivatives of $w$. For the transverse control,
we start from \beq\label{app13} \lambda_{1,z} + m_{B,\Omega}^N \geq
\lambda_{1,z}
 + \frac 12
 \int_{\mathbb R^2_{x,y}\times]-\frac{NT}{2},\frac{N}{2}} |\nabla_{x,y} w -i
\Omega r_\perp w|^2 dx dy\,,
\eeq
which leads to
\beq\label{app14}
 ||  |\nabla_{x,y} w -i
\Omega r_\perp w|\,||^2  \leq 2 m_{B,\Omega}^N\,.
\eeq
For the longitudinal control, we write, for any $\alpha \in [0,1]$
\beq\label{app15}
\lambda_{1,z} + m_{B,\Omega}^N \geq \lambda_{1,z} ||\Pi_N\Psi||^2
 + \frac \alpha 2 || \pa_z w||^2 + \lambda_{N+1,z} (1-\alpha) ||w||^2\,.
\eeq
We determine $\alpha$ by writing
$$
\lambda_{N+1,z} (1-\alpha) =\lambda_{1,z}\,,
$$
hence
\beq\label{app16}
\alpha = 1 - \frac{\lambda_{1,z}}{\lambda_{N+1,z}}\,.
\eeq
So we have
\beq \label{app17}
|| \pa_z w||^2 \leq \frac{2}{\alpha}  m_{B,\Omega}^N \leq
 2 \frac{\lambda_{N+1,z}}{\delta_{N,z}}  m_{B,\Omega}^N\,.
\eeq In the semi-classical regime where we are, this leads to
 the existence of a constant $C$ such that
\beq\label{app18}
|| \pa_z w||^2 \leq C  m_{B,\Omega}^N\,.
\eeq
Using in addition the diamagnetic inequality, we obtain
\beq\label{app19}
|| \nabla |w| ||_2^2 \leq C  m_{B,\Omega}^N\,.
\eeq

As in the other case, we obtain from Sobolev's Inequality the control
 of $w$  in
$L^6$ norm
\beq\label{app20}
|| w||_{6} \leq C ( m_{B,\Omega}^N)^{\frac 12} (1 +
\frac{1}{\delta_z^N})^{\frac 13} \leq \widetilde C ( m_{B,\Omega}^N)^{\frac 12}\,,
\eeq
where we have used that $\delta_z^N >>1$ in the semi-classical
 regime.\\
Using H\"older's inequality, we obtain \beq \label{app21} || w||_4
\leq C ( m_{B,\Omega}^N)^{\frac 12} (\delta_{z}^N)^{-\frac 18}\,.
\eeq We now have all the estimates needed  to mimic the proof of
case
 A.\\

We start from
\beq \label{app22}
\Epsilon(\Psi) \geq \lambda_{1,z} + \Epsilon_B (\Pi_N \Psi)
 - 4 g \int |\Pi_N \Psi|^3 |w|\;dx dy dz\,.
\eeq
We have now to control the third term in \eqref{app22} by the second
term. This is done like in case A in the following way~:
\beq \label{app23}
\begin{array}{ll}
 4 g \int |\Pi_N \Psi|^3 |w|\;dx dy dz &\leq 4 g ||\Pi_N
 \Psi||_4^{3}\; || w||_4\\
& \leq C_1 g^{\frac 14} (\delta_{z}^N)^{-\frac {1}{8}} \left( \Epsilon_B (\Pi_N
 \Psi)\right)^\frac 34 \; ( m_{B,\Omega}^N)^{\frac 1 2}\,.
\end{array}
\eeq
We now use
\beq\label{app24}
\Epsilon_B (\Pi_N\Psi) \geq m_{B,\Omega}^N ||\Pi_N\Psi||_2^4\,,
\eeq
 which together with \eqref{app11} leads to
 \beq\label{app25}
 m_{B,\Omega}^N \leq C (1 + \frac{m_{B,\Omega}^N}{\delta_z^N})
\Epsilon_B (\Pi_N\Psi)\,.
\eeq
This  leads to
\beq\label{app26}
 4 g \int |\Pi_N \Psi|^3 |w|\;dx dy dz
 \leq C_2  g^{\frac 14}( m_{B,\Omega}^N)^{\frac 1 4}  (\delta_{z}^N)^{-\frac {1}{8}}
(1 + \frac{m_{B,\Omega}^N}{\delta_z^N})\Epsilon_B (\Pi_N
 \Psi)\,.
\eeq Using this control, \eqref{app11},  \eqref{app22} and
\eqref{app24}, we have obtained
 the detailed proof of \eqref{minB} in the general case.

\end{proof}
\subsection{On the minimizers of $\Epsilon_B$.}
The next proposition is rather standard and  refers to the case $N=1$ but has its own
interest.
 As a corollary, this will yield an upperbound for $m_{B,\Omega}$.

\begin{proposition}\label{Prop6.10}~\\
The minimizer of $\Epsilon_B$ over the normalized $\psi$'s is unique
 (up to a multiplicative constant of modulus $1$)
and  radial.
\end{proposition}
\begin{proof}~\\
 We first observe that if $\psi$ is a
minimizer
 then $|\psi|$ is also a minimizer. Consequently, we start considering
a non negative minimizer.\\
Now $|\psi|$  is solution of the corresponding Euler equation
 and by the Maximum Principle, $|\psi|$ cannot have a local minimum.
Hence $\psi$ cannot vanish and we can write
$$
\psi = |\psi|\, e^{i\alpha}\,.
$$
Comparing  $\Epsilon_B(\psi)$ and $\Epsilon_B(|\psi|)$ we get
$$
|\nabla \psi| = |\nabla |\psi|\,|\;a.e.
$$
and this implies that $\alpha$ is constant.\\
So we can now assume that $\psi$ is a real positive minimizer. The
Euler-Lagrange equation reads \beq -\Delta \psi + \omp^2 r^2 \psi +
g |\psi|^2 \psi =\lambda \psi \eeq
for some real Lagrange multiplier $\lambda$.\\
Let $\phi$ another positive solution  (with $||\phi||_{L^2}=1$)
 of the Euler-Lagrange equation for a possibly different  $\mu\in \mathbb R$~:
\beq -\Delta \phi + \omp^2 r^2 \phi + g |\phi|^2 \phi =\mu  \phi\,.
\eeq Possibly exchanging the roles of $\psi$ and $\phi$, we can
w.l.o.g assume
 that
$$
\lambda \geq \mu\,.
$$
Let us consider, for some  $\alpha >0$ to be determined,  the
rescaling
$$
\phi (x,y) =\sqrt \alpha \, u (\sqrt{\alpha} x,\sqrt{\alpha} y)\,,$$
we get for $u$ the equation
$$
-\Delta u  + \frac{\omp^2}{\alpha^2}
 r^2 u + g |u|^2 u = \frac \mu \alpha  u \,.
$$
We now choose $\alpha = \frac \mu \lambda$ which leads to \beq
-\Delta u  +\omp^2 \left(\frac{\lambda^2}{\mu^2} \right)
 r^2 u + g |u|^2 u = \lambda   u \,.
\eeq We can now compare $u$ and $\psi$. Let us introduce
 $v =\frac u \psi$ which is a solution of
$$
- \div \left( \psi^2 \nabla v\right) + gv \psi^4 (|v|^2 -1) =
 \omp^2 r^2 v \psi^2 (1 -\frac{\lambda^2}{\mu^2})\,.
$$
Multiplying by $(v-1)_+$  and integrating\footnote{In full rigor,
 we should consider a sequence  $\chi_n (v-1)_+$
 where $(\chi_n)_n$ is a suitable sequence of  cut-off functions and take the limit (See
 \cite{BrOs}).}
we obtain
$$
\int_{\mathbb R^2} \psi^2 |\nabla (v-1)_+|^2 + gv \psi^4
(v-1)_+^2(v+1) \; dxdy
 = \int_{\mathbb R^2} \omp^2 r^2 v \psi^2 (1-\frac{ \lambda^2}{\mu^2})\;
 dx dy\,.
$$
With our assumption on $(\lambda,\mu)$, this implies the vanishing
of
  $(v-1)_+$
 almost everywhere hence $v\leq 1$.\\
This can be reinterpreted as $u\leq \psi$ and the $L^2$
normalization implies $u=\psi$. The lemma is proved.

Finally, one can construct a radial positive solution (by minimizing
 $\Epsilon_B$ over the radial functions). This gives a solution $\phi$
 of the Euler equation which is also  strictly positive. \end{proof}

Observing that, if $\psi$ is radial, we have that
$ \Epsilon_{B,\Omega}(\psi)= \Epsilon_B(\psi)$, this proposition has the following interesting corollary.

\begin{corollary} ~\\
We always have \beq \label{comprot} \inf \Epsilon_{B,\Omega}:=
m_{B,\Omega}  \leq m_B\,. \eeq
\end{corollary}

\subsection{Lower bounds in the TF case ($N\geq 1$)}
We start from a minimizer $(\psi_{\ell,\perp})_\ell$.
Due to the
 normalization, there exists at least one $j$ such that
$$||\psi_{j,\perp}||\geq \frac{1}{\sqrt{N}}$$
Then we write (neglecting the kinetic part)
$$
m_{B,\Omega}^N
 \geq \frac 12(\omega^2-\Omega^2) \int r^2 |\psi_{j,\perp}|^2  + g \int
_{-\frac{NT}{2}}^{\frac{NT}{2}}\int_{\mathbb R^2_{x,y}}
 \left( \sum_{j=0}^{N-1} \psi_j^N(z)\psi_{j,\perp} (x,y) \right)^4\; dz
dx dy\,.
$$
When expanding $ \left( \sum_{j=0}^{N-1} \psi_j^N(z)\psi_{j,\perp}
(x,y) \right)^4 $, the mixed terms
 are exponentially small (see Subsection \ref{sectionDNLS}) in comparison to
 $\sum_j ||\psi_{j,\perp}||_{L^4}^4$, hence we get, for some $\alpha
>0$,
$$
m_{B,\Omega}^N  \geq \frac12 (\omega^2-\Omega^2) \int r^2 |\psi_{j,\perp}|^2
 + g \int_{-\frac{NT}{2}}^{\frac{NT}{2}} \psi_0^N(z)^4 dz
(\int (\psi_{j,\perp})^4 dxdy)\; (1- \exp - \frac{\alpha}{\epsilon})\,.
$$
We now use \eqref{quartica}, to obtain
$$
\begin{array}{ll}
m_{B,\Omega}^N  &\geq\frac 12 (\omega^2-\Omega^2) \int r^2 |\psi_{j,\perp}|^2
 + g  \int_{-\frac{T}{2}}^{\frac{T}{2}} \phi_1(z)^4 dz
(\int \psi_{j,\perp}^4 dxdy) (1- \exp - \frac{\alpha}{\epsilon})\\&
= \frac 12(\omega^2-\Omega^2)\int r^2 |\psi_{j,\perp}|^2
 + \widetilde g
(\int \psi_{j,\perp}^4 dxdy) (1- \exp - \frac{\alpha}{\epsilon})\\&
\geq  \left(\frac 12(\omega^2-\Omega^2)\int r^2 |\psi_{j,\perp}|^2
 + \widetilde g
(\int \psi_{j,\perp}^4 dxdy) \right) (1- \exp -
 \frac{\alpha}{\epsilon})\\&
\geq \frac{1}{N^2} (1- \exp -
 \frac{\alpha}{\epsilon})  \inf_{\psi,||\psi||=1} \left(\frac 12 (\omega^2-\Omega^2)\int r^2 |\psi|^2
 + \widetilde g
(\int \psi^4 dxdy) \right)
 \,.
\end{array}
$$
One can then use the asymptotics
 obtained in the proof of \eqref{asympcasBOmegaN=1}
 to get, under Assumption \eqref{controleOmega},  the existence of
 $C_{N,\beta}>0$ such that, as $\ep$ tends to \ 0 and
$ \widetilde g$ to $\infty$, \beq\label{lbmbntf} m_{B,\Omega}^N
 \geq \frac{1}{C_{N,\beta}} \omp \sqrt{\tilde g}\,.
\eeq

\section{Tunneling effects for the non-linear models}\label{stunneling}
This is only in this section that we will exhibit the role of these
localized $(NT)$-periodic Wannier functions.

\subsection{Towards the DNLS model.}\label{sectionDNLS}
\subsubsection{Preliminaries}
We have already proved for any $N\geq 1$ the rough estimates on $m_A^N$
 allowing to justify the longitudinal model. We have established or
 announced better asymptotics at the price of
 stronger assumptions.\\
 Our aim in this
 section is to discuss possible asymptotics for $m_A^N$
 in the case when $N>1$, which will involve the tunneling effect. Although we have no
  final result
 on this part, we would like to prove how we reach a
 familiar
 model considered by the physicists~: the DNLS model.\\
In particular we will describe in Proposition \ref{justifsnoek}
 under which assumptions one can get
 a simplified model.

 We consider on $]-\frac {NT}{2},\frac{NT}{2}[$ the
$(NT)$-periodic
  problem for  the operator
 $-\frac{d^2}{dz^2} + W_\ep(z)$
with $W_\ep(z) = \frac{{\bf w}(z)}{\ep^2} $, $w$ satisfying Assumption
\eqref{assw}. Here  we always work in the semi-classical regime.\\

The starting point in this subsection is that we replace the issue
 of minimizing $\Epsilon_A^{N,\ep,\widehat g}$ on the $(NT)$-periodic
 $L^2$-normalized functions by restricting the approximation to the
 eigenspace ${\rm Im}\, \pi_N$  associated with the first $N$  eigenvalues
of the linear
 problem.

\subsubsection{Projecting on the eigenspace ${\rm Im}\, \pi_N$}
Our aim is to analyze the reduced functional \beq \mathbb C^N\ni \cb
=(c)_{j=0,\dots,N-1}\mapsto  \Epsilon_A^{N,\ep,\widehat
  g,red}(\cb)=
\Epsilon_A^{N,\ep,\widehat g} (\sum_{j=0}^{N-1} c_j \psi_j^N)\,, \eeq
 where  $\Epsilon_A^{N,\ep,\widehat g}$ is in fact $\Epsilon_A^{N}$
 given
 in \eqref{EpA} with the explicit notation of
 the dependence of the parameters and the $\psi_j^N$ are the $(NT)$-periodic Wannier
 functions. When $N=1$, the error which is done has been
 estimated
 in \eqref{mAimp}  under the assumption
 that
 $\widehat g \epsilon^\frac 12$ is small,
 i.e. \eqref{condbizb}. Replacing
 in the argument the projection on the first eigenspace by $\pi_N$,
 the same result holds true for $N>1$.

We now concentrate our discussion to the model obtained after this
first approximation.
More specifically we are interested in the asymptotics of the infimum 
of this functional.\\

Here natural approximations of this reduced functional appear. Each of
 these approximations gives a corresponding  approximation of the
 infimum of the reduced functional, which is defined by~:
\beq m_A^{N,(0)}(\ep, \widehat g) := \inf_{\{\cb\,|\,\sum_{j=0}^{N-1}
|c_j|^2 =1 \} }  \Epsilon_A^{N,\ep,\widehat
  g,red}(\cb)\,.
\eeq Let $\lambda_{1,z}^{N}=\lambda_{1,z}$
 be the bottom of the $(NT)$-periodic spectrum of $H_z$ on
 $[-\frac{NT}{2}, \frac{NT}{2}]$ (with $N$ minima). So strictly
 speaking, we can start the analysis of this first approximate model
 only under Condition \eqref{condbizb}.
\begin{proposition}~\\
Under condition \eqref{condbizb}
\beq
m_A^N (\epsilon,\widehat g)= m_A^{N,(0)} (\epsilon,\widehat g)
 + \mathcal O (\widehat g^\frac 32 \ep^{-\frac 14})\,.
\eeq
\end{proposition}
 One can nevertheless imagine
 that the information obtained in the next subsubsection
 is valid in a more general context (maybe by choosing other localized
 Wannier functions).
 We now analyze various approximations of
 $m_A^{N,(0)}(\ep, \widehat g)$.

\subsubsection{Neglecting the tunneling}
 Neglecting the tunneling effect,  we are lead to  the minimum of the functional
$\Epsilon_A^{N,\ep,\widehat g,(1)}$
$$
\mathbb C^N\ni\cb \mapsto \Epsilon_A^{N,\ep,\widehat g,(1)}(\cb):=
\lambda_{1,z} \left(\sum_{j=0}^{N-1} |c_j|^2\right) + \widehat g \,
(\sum_{j=0}^{N-1} |c_j|^4)\; (\int_{-\frac{NT}{2}}^{\frac{NT}{2}}
|\psi_0^N(z)|^4 \, dz)\,,
$$
over the $c$'s such that
$$
\sum_{j=0}^{N-1} |c_j|^2=1\,.
$$
Observing (see \cite{DiSj}), that 
\beq\label{quartica}
\int_{-\frac{NT}{2}}^{\frac{NT}{2}} |\psi_0^N(z)|^4 \, dz =
\int_{-\frac T2}^{\frac T2}\phi_1(z)^4\, dz  + \widetilde \Og(\exp -
\frac {S}{2 \ep})\,, \eeq
 where $\phi_1$ is the groundstate of
the $T$-periodic problem, the minimum of this  approximate functional,
which is attained for $c_j = N^{-\frac 12}$,  is 
\beq m_A^{N,(1)} =
\lambda_{1,z}+ \frac{\widehat g}{N}
\int_{-\frac{NT}{2}}^{\frac{NT}{2}} |\psi_0^N(z)|^4 \, dz\,. \eeq

So as a first approximation, we have obtained
\begin{proposition}~\\
$$
m_A^{N,(0)}(\ep, \widehat g) = \lambda_{1,z} + \frac{ \widehat g}{N}
(\int_{-\frac{NT}{2}}^{\frac{N T}{2}} |\psi_0^N(z)|^4 \, dz) +
(\widehat  g +1) \; \widetilde\Og (\exp - \frac {S}{\ep})\,,
$$
or
\beq\label{asympassezprecise}
m_A^{N,(0)} (\ep,\widehat g)= \lambda_{1,z} +  \frac{\widehat
g}{N}(\int_{-\frac T2}^{\frac T2} \phi_1(z)^4\,dz)
 + \widehat  g \; \widetilde\Og (\exp - \frac {S}{2 \ep}) + \widetilde\Og (\exp - \frac {S}{\ep})
\,.
\eeq
\end{proposition}
The definition of $\widetilde \Og$ is given in \eqref{defOtilde}. If
we apply this result to our context with $\widehat g =\omp g$, this
yields  information on the behavior of $m_A^{N,(0)}$  independently
of Assumption \eqref{condbizb}.

\subsubsection{Taking into account the tunneling}
If we keep the main tunneling term, we get the following more
accurate approximating functional
\begin{equation}
\begin{array}{l}
\mathbb C^N\ni \cb \mapsto
\Epsilon_A^{N,\ep,\widehat g,(2)}(c)\\
\quad :=

\widehat \lambda_1 \left(\sum_{j=0}^{N-1} |c_j|^2\right) -  \tau
\Re\left  (\sum_{j=0}^{N-1} c_j\,\overline{c_{j+1}}\right) +
\widehat g \,  (\sum_{j=0}^{N-1} |c_j|^4)\;
(\int_{-\frac{NT}{2}}^{\frac{NT}{2}} |\psi_0^N(z)|^4 \, dz)\,.
\end{array}
\end{equation}
Here $\tau$ is the hopping amplitude introduced around
\eqref{asympoft},  $\widehat \lambda_1$
is
 the lowest eigenvalue corresponding to the Floquet condition $k=\frac N2$
 for the linear problem on $]-\frac T2, \frac T2[$, which is
exponentially closed to $\lambda_1$ 
 and we take  the convention that $c_{N}=c_0$.\\
The quadratic form corresponds to the approximation in the first
band~: \beq \mathbb C^N\ni \cb \mapsto \widehat \lambda_1
\left(\sum_{j=0}^{N-1} |c_j|^2\right) - \tau  \Re\left
(\sum_{j=0}^{N-1} c_j\,\overline{c_{j+1}}\right) \eeq which can be
shown to be correct modulo $\widetilde \Og (\exp - \frac{2S}{\ep})$.
\begin{remark}~\\
This time the minimizer could depend on $\widehat g$ !! This is the
kind of problem
 which is analyzed in \cite{KMPS} and in Subsection \ref{lastsection}.
\end{remark}

\paragraph{Discussion about the justification of $\Epsilon_A^{N,\ep,\widehat g,(2)}$}~\\
One can wonder why we
 forget some terms in the computation. Let us do this more
 carefully.
To be consistent with what we forget in the linear case (terms of
order
 $\Og (\tau ^2)$), we show first that one can  approximate\footnote{We use here the
 assumption that the potential and hence $\psi_0^N$ is even. We recall
 also that the $\psi_j$ are real.}
 $\left( \int_{-\frac{NT}{2}}^{\frac{NT}{2}} |\sum_{j=0}^{N-1} c_j \psi_j^N(z)|^4 \, dz\right)$ by
\begin{equation}\label{Err1}
\begin{array}{l}
\left( \int_{-\frac{NT}{2}}^{\frac{NT}{2}} |\sum_{j=0}^{N-1} c_j \psi_j^N(z)|^4 \, dz\right)
 = \\
\quad (\sum_{j=0}^{N-1} |c_j|^4) (\int_{-\frac{NT}{2}}^{\frac{NT}{2}} |\psi_0^N|^4 dz)\\
 + \sum_{j=0}^{N-1}\left( (|c_j|^2 +|c_{j+1}|^2) (c_j \,
\overline{c_{(j+1)} }+c_{j+1}\,\overline{c_{(j)} }
 ) \;(\int_{-\frac{NT}{2}}^{\frac{NT}{2}} \psi_0^N(z)
 |\psi_0^N(z)|^2\cdot \psi_1^N(z) dz)\right) \\
 + \widetilde \Og (\tau^2)\,.
\end{array}
\end{equation}
This first approximation is based on the following lemma.
\begin{lemma}~\\
 $$
\int_{-\frac{NT}{2}}^{\frac{NT}{2}}\psi_0^N (z)^2 \psi_1^N(z)^2 dz =  \widetilde \Og (\exp - \frac
        {2S} \ep)\,.
$$
\end{lemma}
This is based on the property that, for all $\eta >0$,  there exists
$C_\eta$ such that
\begin{equation}\label{AgmDec}
|\psi^N_0(z)|\leq C_\eta \exp \frac \eta h \;  \exp - \frac  1 \ep
d^{mod}_{Ag}(z)\,,
\end{equation}
where $d^{mod}_{Ag}(z)$ is an even function
 such that
$$
d^{mod}_{Ag}(z,0) = 2  \int_0^z \sqrt{{\bf w}(t)} \,dt\,,\; \mbox{ for }
z\in [0, T[\,,
$$
and such that $
d^{mod}_{Ag}(z,0)$ is increasing for $z\geq 0$.\\
On the contrary, this  is a priori  unclear\footnote{In \cite{KMPS},
p.~5, between
 formulas (18) and  (19), the term $\widehat \tau$ is discussed; see
 also p.~6 around formula (20).}   why one could forget
 terms like
\begin{equation}\label{defthat}
  \widehat \tau = \widehat g \int_{-\frac{NT}{2}}^{\frac{NT}{2}} \psi_0^N(z)^3 \psi_1^N(z) dz \,.
\end{equation}
(where  we recall that ${\bf w}$ is even  by Assumption \eqref{assw} and
that
 this
 implies  $\psi_0^N$ even and real).
This term is a priori of the same order
 as $\tau $. We have indeed
\begin{lemma}\label{lemmixte}~\\
 \beq\label{termemixte}
\int_{-\frac {NT}2}^{+\frac {NT}2} \psi^N_0(z)^3 \psi^N_1(z)\; dz
 = \widetilde \Og (\exp - \frac {S} \ep)\,.
\eeq
\end{lemma}
 Due to the decay estimates \eqref{AgmDec} for these $(NT)$- Wannier functions, the term to
 integrate in \eqref{termemixte}  decays like
 $$
\widetilde \Og \left(\exp -\frac 1 \ep  \left(3  d_{Ag}^{mod}(z) +
d_{Ag}^{mod} (z- T)\right) \right) \,,
$$
so the main contribution comes from the origin
 and has the same size as $\exp - \frac S \ep$.\\

So it is necessary to be careful\footnote{We thank
 M.~Snoek for kindly answering our questions on this problem.}, if one wants to
 neglect $\widehat \tau$.\\

Let us now try to estimate $\int_{-\frac {NT}{2}}^{+\frac {NT} 2}
 \psi^N_0(z)^3 \psi^N_1(z)\; dz$ as $\ep
\rightarrow 0$ more precisely. Heuristically, one can try to use a
WKB
 approximation, this is available for $\psi_0^N$ in the neighborhood
 of $0$ but
unfortunately, we do not have a good WKB approximation of $
\psi^N_1(z)$ close to the origin, as observed in Subsection
\ref{sstunneling} (see \eqref{constance}). So we have no obvious
main term
 for the asymptotic behavior of  $\int_{-\frac {NT}{2}}^{+\frac {NT} 2}
 \psi^N_0(z)^3 \psi^N_1(z)\; dz$. A reasonable guess (which is
 implicitly used
 by the physicists) should be to
 suggest the following conjecture.
\begin{conjecture}~\\
\begin{equation}\label{asymptofthat}
\widehat \tau = \widehat g \; o( \tau )\,,
\end{equation}
as $\ep \ar 0$.\end{conjecture}

The weaker  mathematical result, which is obtained from Lemma
\ref{lemmixte} by the considerations above  using
Helffer-Sj\"ostrand techniques \cite{DiSj},  is the following
proposition.
\begin{proposition}\label{justifsnoek}~\\
Under the assumption that there exists $\eta >0$ such that,
\begin{equation}
0\leq \widehat g  \exp \frac \eta \ep \leq 1 \,,
\end{equation}
   then
\beq m_A^{N,(0)} = m_A^{N,(2)} + o(\tau )\,. \eeq holds.
\end{proposition}

This gives a motivation for the analysis of
 the DNLS model of \cite{STKB} (with an extra term in $\lambda
 \sum_{j=0}^{N-1}
 |c_j|^2$).\\
If we consider the $(NT)$-periodic Floquet problem, we arrive
 naturally to questions analyzed in  \cite{KMPS}
(16-17-18), and the remark after (21) in this paper.

\subsection{On approximate models in  case $B$ : towards Snoek's model}\label{SubsectionSnoek}

Using the basis of the ($NT$)-Wannier
 approach), we can  consider $\mathcal E_B^{N}$
 introduced in \eqref{plussimple} and consider the decomposition
$$
\mathcal E_B^{N}(\psi_{0,\perp},\cdots,\psi_{N-1,\perp}) : =
\mathcal E_B^{N\,'}(\psi_{0,\perp},\cdots,\psi_{N-1,\perp})
 + g ||\sum_{j=0}^{N-1}\psi_j^N(z) \psi_{j,\perp} (x,y)||_{L^4}^4\,.
$$

We now use various approximations related to the analysis of the
$z$-problem ($(NT)$-Wannier functions). We get

$$
\begin{array}{l}
\mathcal E_B^{N\,'}(\psi_{0,\perp},\cdots,\psi_{N-1,\perp})
\\
 \quad \sim s \sum_{j=0}^{N-1} ||\psi_{j,\perp}||^2 +
  t \sum_{j=0}^{N-1} \left(\langle \psi_{j,\perp},\psi_{j+1,\perp}\rangle
 + \langle \psi_{j,\perp},\psi_{j-1,\perp}\rangle\right)\,,
\end{array}
$$
and
$$
\begin{array}{l}
 g ||\sum_{j=0}^{N-1} \psi_j(z) \psi_{j,\perp} (x,y)||_{L^4}^4
 \sim g ||\psi_0||_{L^4}^4 \sum_{j=0}^{N-1} ||\psi_{j,\perp}||_{L^4}^4\,.
\end{array}
$$
So the approximate functional becomes

\begin{equation}\label{EpsilonB}
\begin{array}{ll}
\Epsilon_{B}^{N,approx} ((\psi_{j,\perp})_j) &=  \sum_{j=0}^{N-1}
\int_{\mathbb R^2} \left( \frac 12 |\nabla
\psi_{\perp, j}|^2 + V(x,y)|\psi_{j,\perp}(x,y)|^2\right) dx dy
\\&\quad  +
 s \sum_{j=0}^{N-1} ||\psi_{j,\perp}||^2\\&\quad+  t
\sum_{j=0}^{N-1}\left(\langle \psi_{j,\perp},\psi_{j+1,\perp}\rangle
 + \langle \psi_{j,\perp},\psi_{j-1,\perp}\rangle\right)\\&
\quad +  \widetilde g \sum_{j=0}^{N-1} ||\psi_{j,\perp}||_{L^4}^4\,,
\end{array}
\end{equation}
which should be minimized over the $(\psi_{j,\perp})_j$
 such that
$$
\sum_{j=0}^{N-1} ||\psi_{j,\perp}||^2 =1\,.
$$
This is the model described by  Snoek \cite{Sn}.

Starting from this model, one can,  depending on the size of the
various parameters, come back in some case to  the situation when
 $(\psi_{j,\perp})_j$ is of the form
 $c_j \psi_\perp$, with $\sum_{j=0}^{N-1} |c_j|^2=1$.  In this case, we come
back to the results of the
previous subsection. In the other cases, the problem seems completely
open.

\subsection{Spatial period-doubling in Bose-Einstein condensates
 in an optical lattice}\label{lastsection}
Here we mainly follow \cite{MNPS}. We look at the discrete model for which
these authors refer to \cite{STKB}.

The Hamiltonian for the discrete model is formally
 defined\footnote{after substraction of a term in the form $\sigma
 \sum_j |c_j|^2$,}
 as
$$
H (\cb) = - \tau  \sum_{j} (\overline{c_j} c_{j+1} +c_j\overline{c_{j+1}} )
 + I \sum_{j}|c_j|^4\,,
$$
where $\cb =(c_j) \in \ell^2(\mathbb Z;\mathbb C)$, $\tau \geq 0$, $I  \geq 0$.\\

This is a particular case of  the so-called $DNLS$ model\footnote{In
  the general case one should add
 the term $\sum_j \ep_j |c_j|^2$
 when the trapping is present the authors propose $\ep_j = \tilde
  \omega j^2$
 with $\tilde \omega$ proportional to $\omega_z$.} (Discrete Non
  Linear Schr\"odinger model).
But we will immediately reduce our analysis to the $(NT)$-periodic
  problem
 and Floquet variants of this problem.\\
We will restrict the sum above to $j=1,\dots, N$, where $N$
 is a fixed positive integer. But  for defining the tunneling term
$ (\overline{c_j} c_{j+1} +c_j\overline{c_{j+1}} ) $,
we use the Floquet condition
\beq \label{FloqNT}
 c_{N +1} = \exp( ik  N) \;  c_{1}\,.
\eeq
\begin{remarks}~
\begin{itemize}
\item
\cite{MNPS} looks at the particular case $N=2^p$.
\item
 When $k=0$, this is the
natural $(NT)$-periodic problem.\item
 Note that it is the standard Floquet problem
 only for $N=1$.\item
Note that we forget the term $\lambda \sum_j |c_j|^2$ which
 corresponds only to
 a shift of the energy.
\end{itemize}
\end{remarks}
Hence we would like to minimize $H^{N,k}(\cb)$
\beq
H^{N,k} (\cb) :=
 - \tau\,  \sum_{j=1}^N (\overline{c_j} c_{j+1} +c_j\overline{c_{j+1}} )
 + I\, \sum_{j=1}^N |c_j|^4\,,
\eeq
 over the  $\cb$'s normalized in $\mathbb C^N $
 by
\begin{equation}
||\cb||^2 = N_c\,,
\end{equation}
and the $k$-Floquet condition \eqref{FloqNT}.\\
Moreover for future use we introduce the strictly positive parameter
\begin{equation}
\nu  = N_c/N >0\,.
\end{equation}
\begin{remark}~\\
In the preceding sections we were taking $N_c=1$. Up to a change of
the
 parameter $\tau  $, one can always reduce to the general case to this
 situation.
 A difference could occur if we take $\nu$ fixed
 and $N\ar +\infty$.
\end{remark}

We will then be interested in the analysis of the energy per particle
\beq \label{deefenerg}
E (\tau , I ,\nu, N,k) = \frac{1}{N_c} \inf_{||\cb||^2 =N_c} H^{N,k}
(\cb)\,.
\eeq

 Writing
\begin{equation}
c_j = \exp i\,k j\; g_j
\end{equation}
then from \eqref{FloqNT}
\begin{equation}\label{perNT}
g_{1+N} = g_1\,.
\end{equation}
The case when $N=1$ corresponds to the usual Floquet condition.\\
Writing that $\cb$ is a critical value of $H^{N,k}$, we get
 that for some $\mu$ (which is called the chemical potential
 or the Lagrange multiplier in mathematics)
\begin{equation}
2 I |c_j|^2 c_j - \tau  \left( c_{j+1} + c_{j-1}\right)= \mu c_j\,,
\end{equation}
for $j=1,\dots, N$, with condition \eqref{FloqNT}.
This becomes in terms of the $g_j$'s
\begin{equation}
2I |g_j|^2 g_j - \tau  \left(\exp ik \; g_{j+1} + \exp(- ik)
\;g_{j-1}\right)= \mu g_j\,,\end{equation}
 for $j=1,\dots, N$, with the $N$-periodic convention \eqref{perNT}.

\paragraph{The case $ N=1$}~\\
In this case, we have simply one equation
\begin{equation}
2I  |g_1|^2 g_1 - 2 \tau  \cos (k) \, g_1= \mu\,  g_1\,,
\end{equation}
 with
\begin{equation}
|g_1|^2 = \nu\,.
\end{equation}
We find immediately that
\begin{equation}
\mu =  -2 \tau  \cos( k) + 2 I  \nu
\end{equation}
and
\begin{equation}
g_1 = \nu \exp i \phi_1 \,.
\end{equation}
The energy per particule $E$
 is then equal to
\begin{equation}
E =\frac{1}{N_c}  H^{1,k} (\cb) =  -2 \tau  \cos( k) + I  \nu\,.
\end{equation}
For this choice of $N$, we have recovered exactly what we have found
in the linear case $I =0$. The effect of the non-linear term just
creates a $k$-independent shift of the energy. Note also that as a
function of $k$, the energy is minimal for $k=0$, which is the
periodic case.\\

\paragraph{The case $ N=2$}~\\
In this case, we get, using the periodicity assumption,
  the following system of equations
\begin{equation}\label{dnls1}
\begin{array}{l}
2I  |g_1|^2 g_1 - 2 \tau  \cos (k)\, g_2 = \mu g_1\,,\\
- 2 \tau  \cos( k) \, g_1  + 2I  |g_2|^2 g_2 =  \mu g_2\,,
\end{array}
\end{equation}
with the normalization condition
\begin{equation}
|g_1|^2 + |g_2|^2 = N_c = 2 \nu\,.
\end{equation}
We write
\begin{equation}
g_j = |g_j|\exp i \varphi_j\,,\; \mbox{ for } j=1,2\,.
\end{equation}
A suitable combination of the two lines
 gives
\begin{equation}\label{condinec}
2 I  \left( |g_1|^2 - |g_2|^2\right)
 = 2 \tau  \cos (k) \left( \frac{|g_2|}{|g_1|} \exp i (\varphi_2-
 \varphi_1)\;
 - \;  \frac{|g_1|}{|g_2|} \exp-  i (\varphi_2-
 \varphi_1)\right)\,.
\end{equation}
If we observe that the left hand side is real, we get
 \begin{equation}
0 = 2 \tau  \cos (k) \left( \frac{|g_2|}{|g_1|}
 + \;  \frac{|g_1|}{|g_2|} \right) \sin  (\varphi_2-
 \varphi_1)\,.
\end{equation}
The real part of \eqref{condinec} is in any case given by
\begin{equation}
I  \left( |g_1|^2 - |g_2|^2\right)
 =  \tau   \cos (k) \left( |g_2|^2 - |g_1|^2\right) \frac{1}{|g_1||g_2|}
 \cos  (\varphi_2- \varphi_1)\,.
\end{equation}

{\bf We meet three cases.}\\

{\bf Case 1}~\\
We first observe that the solutions corresponding to the
 case $N=1$ are recovered by taking $\varphi_1=\varphi_2$
 and $|g_1|=|g_2|$.\\
Another family of solutions is obtained by taking $\varphi_1=\varphi_2
+\pi$
 and $|g_1|=|g_2|$. This corresponds to an ``antiperiodic''
 solution over two periods.

The solutions such that $|g_1| =|g_2|=\sqrt{\nu}$ seem to be simply
deformations
 of the case $I=0$.\\
 At least for $I$ small, this is indeed
 also a consequence of the implicit  function theorem
 when $\cos k \neq 0$.\\
The energy per particle $E$ is the same as for $N=1$.

{\bf Case 2}~\\
If we  assume that
\begin{equation}
|g_1|\neq |g_2|\,,
\end{equation}
and
\begin{equation}
\tau  \cos k \neq 0\,,
\end{equation}
then the previous necessary conditions become
 first
$$
\sin (\varphi_2-\varphi_1) =0\,,
$$
hence
$$
\varphi_2=\varphi_1 \;mod (\mathbb Z \pi)\,,
$$
and secondly (using the first one)
\begin{equation}
2I  |g_1||g_2|
 =\pm  2 \tau   \cos k \,.
\end{equation}
Using the normalization of $g_1$ and $g_2$, we get
 as a necessary condition
\begin{equation}
  |\tau  \cos k |\leq |I|\nu\,.
\end{equation}

If these conditions are satisfied, we can find in function of the
 sign\footnote{It seems that in our problem we have $I\geq 0$ and
 $\tau >0$.} of $\cos k$, a unique
 pair $g_1$, $g_2$ (up to a multiplication by $e^{i\theta}$).  Coming back to the initial system
 of equations leads to the determination of $\mu$
 which is given
 by
\begin{equation}
\mu = 4 I  \nu\,,
\end{equation}
and of the energy  $E$
 which is given by~:
\begin{equation}
E =\frac{\tau^2}{I\nu}\cos^2k\,+\,  2I   \nu\,.
\end{equation}

{\bf Case 3}~\\
It remains to consider the degenerate situation. If
\begin{equation}
\cos k =0\,,
\end{equation}
i.e. if
\begin{equation}
|k|= \frac{\pi}{2}\,,
\end{equation}
Then
\begin{equation}
g_1 = \nu \exp i \varphi_1 \,,\; g_2 = \nu  \exp i \varphi_2\,,
\end{equation}
is a solution for any pair  $(\varphi_1, \varphi_2)$. The
corresponding $\mu$
 is given by
\begin{equation}
\mu = 2I \nu \,,
\end{equation}
and
\begin{equation}
E =  I \nu \,.
\end{equation}
This is the same energy that the one found for the usual Bloch state
(with the same $k$)
 but note that the $g_j$ are no more $1$-periodic ($g_{j+1} \neq g_j$).
\begin{question}~\\
It is unclear in the discussion what is the status of $k$. Are we
interested in minimizing over $k$? But in this case $k=0$
 and $N=1$
 seems optimal in the sense that they give the lowest $E$.
\end{question}

\begin{remark}~\\
An interesting problem, which is discussed in \cite{MNPS}, is the analysis of the stability, looking at the
linearized corresponding problem.
\end{remark}

\begin{remark}~\\
Of course the analysis of more general $N$'s would be quite
interesting.
 A few numerical results are given in \cite{MNPS} corresponding to $N=4$.
\end{remark}


\section{Other   optical lattices functionals}\label{s4}
In this section, we discuss the choice of analyzing periodic
boundary conditions in the $z$ direction and the possibility of
stating the problem differently. We compare the $(NT)$-periodic
problem to the $T$-periodic problem and discuss shortly  the question $N\ar
+\infty$.\\
\subsection{Summary of the linear case}\label{sspc}
 We summarize what we have obtained in the linear
situation. Different techniques can be used for determining the
bottom of the spectrum of $H_z$ and then of $H$, but the ground
state energies always coincide.
\begin{enumerate}
\item Minimize the functional
\begin{equation}\label{defQ}
\psi \mapsto Q(\psi):=\langle H_z \psi\,|\,\psi
  \rangle_{L^2(\mathbb R)}\,,
\end{equation}
 over $L^2$-normalized  $\psi$'s in $ C_0^\infty (\mathbb R)$ (or
  in $H^1(\mathbb R)$). In
  this case, the minimization gives the ground state energy but there
  is no minimizer in the form domain of the operator !
\item Minimize the functional
\begin{equation} \label{defQper}
\psi \mapsto Q^{per}(\psi) =\int_{-\frac T2}^{+\frac
T2}\left(\frac{1}{2}
 |\psi'(z)|^2
 + W_\ep (z)| \psi(z)|^2\right)\, dz \,,
\end{equation} over the
 $L^2$ normalized  $C^\infty$ $T$-periodic functions $\psi$'s
(or on $H^{1,per}$)~. Here we integrate over one period~! The
minimization will give
 the ground state energy of the periodic operator and the minimizer of
 the functional will be the ground state.
\item Minimize the functional
\begin{equation}\label{defQNper}
\psi \mapsto Q^{per, N}(\psi) =\int_{-\frac{NT}{2}}^{\frac{NT}{2}}
\left(\frac{1}{2 }
 |\psi'(z)|^2
 + W_\ep (z)| \psi(z)|^2\right)\, dz \,,
\end{equation}
 over the $L^2$-normalized $C^\infty$ $(NT)$-periodic functions (or on $H^{1,N,per}$).  Here we integrate over $N$ periods~! The
minimization will give
 the ground state energy of the periodic operator and the minimizer of
 the functional will be again  the $T$-periodic ground state.
\item Minimize over $k\in [0,2\pi/T[$, the infimum of the functional
\begin{equation}\label{defQFloq}
Q^{Floq,k}(\psi) =\int_{-\frac T2}^{+\frac T2} \left(\frac{1}{2}
 |\psi'(z)-i k \psi(z)|^2
 + W_\ep (z) |\psi(z)|^2\right)\, dz
\end{equation}
 over the $L^2$-normalized $C^\infty$ $T$-periodic functions $\psi$ (or on
 $H^{1,per}).$
 Here we integrate over one period~! The
minimization will give
 the ground state energy of the periodic operator (i.e. $k=0$) and the minimizer of
 the functional will be the periodic ground state (corresponding to $k=0$).

\item
Minimize over the space generated by the Wannier
  functions (identified with $\ell^2(\Gamma,\mathbb C)$). This
  corresponds to the reduction to the first band and
  leads to the analysis of
\begin{equation}\label{defQWan}
\ell^2(\Gamma,\mathbb C)\ni \cb= (c_\ell)_{\ell \in \Gamma }
 \mapsto \sum_{\ell,m} \widehat \lambda_1(\ell -m) c_\ell
 \overline{c_m}\,,
\end{equation}
with $\sum_\ell |c_\ell|^2=1$.\\

\item
Minimize over the space generated by the $(NT)$-Wannier
  functions (identified with $\ell^2(\Gamma^N,\mathbb
  C)=\mathbb C^N$). This
  corresponds to the reduction to the spectral space
 attached to the first $N$  eigenfunctions of the $(NT)$-periodic problem
 living in the first band.
\end{enumerate}
{\bf We emphasize, that, in each case, we get the same ground state energy.}\\
As already mentioned, the $3$-dimensional linear case introduced in
\eqref{3.1} is easily reduced to the one-dimensional case $H_z$.\\
\begin{remark}~\\
The reader could be astonished that we discuss only the case
$\omega_z =0$. This is a current assumption in the physical
literature. Mathematically, there is a dramatic change in the nature
of the spectrum. The problem could become quite difficult in some
regimes
 of $\omega_z$ and this will not discussed in this paper. Let us
 nethertheless make a few comments.When $\omega_z\neq 0$, the spectrum
 of $H_z$ on the line becomes
 indeed discrete. By monotonicity, the bottom of the spectrum is above
 $\lambda_{1,z}$. One can also get  an upper bound
 by computing the energy of a suitable quasimode (or more simply of
 the ground state
 of the linear problem) but this can only be
 good
 in  some asymptotical regime.
\end{remark}


\subsection{The $3D$- functionals}\label{ss3D}
 Let us consider the fully
non-linear problems and try to implement some of the results
obtained for the linear models. We take $\omega_z=0$,
 we have to analyze (see \eqref{energyground}) the functional
\begin{multline}\label{defQBE}
\Psi \mapsto Q_{\Omega}(\Psi):=\\ \int_{\mathbb R^3}
   \left(\frac{1}{2} |(\nabla_{x,y,z}-i\Omega\times  r) \Psi({\bf r})
   |^2
 +( V({\bf
   r})+ W_\epsilon (z)  -\frac 12 \Omega^2 r^2)|\Psi({\bf r})|^2
   +  g |\Psi({\bf r})|^4 \right) \;dx dy dz \,.
\end{multline}
We denote by $ \mathcal D_{\Omega}$ the natural maximal form domain
 of the form $Q_{\Omega}(\Psi)$ in $L^2(\mathbb R^3)$,
 that is
\begin{equation}\label{natdomr3}
 \mathcal D_{\Omega}=\{u\in H^1(\mathbb R^3)\,,\, xu \in L^2(\mathbb R^3)\,,\, y
 u\in L^2(\mathbb R^3)\}\,.
\end{equation}
We denote  the intersection of the $L^2(\mathbb
R^3)$-unit ball with $ \mathcal D_{\Omega}$~:
$$
 \mathcal S_{\Omega} =\{\Psi \in  \mathcal D_{\Omega}\,,\;
 ||\Psi||_{L^2} =1\}\,.
$$
 We call the infimum of this functional
\begin{equation}
E_\Omega := \inf_{\Psi\in  \mathcal S_{\Omega} }Q_{\Omega}(\Psi)\,.
\end{equation}
Because there is no harmonic trapping in the $z$-variable
  the choice of the
 condition of normalization in $L^2(\mathbb R^3)$  is questionable.\\

The $T$-periodicity of the  optical lattice in the $z$-variable
suggests
  to consider other functionals, where we
integrate
 over $\mathbb R_{x,y}^2 \times ]-\frac T2, \frac T2[$ (or more
 intrinsically over $\mathbb R_{x,y}^2 \times (\mathbb R/T \mathbb
 Z)\,$), or over  $\mathbb R_{x,y}^2 \times ]- \frac {NT}{2},  \frac {NT}{2}[$
 for some integer  $N$, and where the variational space has to be defined
 suitably (periodic conditions or Floquet conditions). We refer to
 Subsection \ref{sspc}, for the discussion done in the
 linear case. In the non-linear case, this has led to the introduction
 of the ''periodic'' Bose-Einstein  functional  (see
 \eqref{defQBEper})
We will denote by $ Q_{\Omega}^{per,N}$ the functional obtained by
integration over $N$ periods (see \eqref{defQBEperN}).\\
We call $\mathcal D_{BE}^{per}$, the natural maximal form domain
 of the form $ Q_{\Omega}^{per}$. It corresponds  to the
 distributions $\Psi$ in $H^1_{loc}(\mathbb R^3)$, satisfying \eqref{propper}
and  such that, the restriction $\Psi_T$ to $ \mathbb
R^2_{x,y}\times
 ]-\frac T2, \frac T2[$,
 satisfies :
\begin{equation}\label{natdom}\Psi_T
\in H^1( \mathbb R^2_{x,y}\times ]-\frac T2, \frac T2[)\,,\;
 \sqrt{x^2+y^2} \;\Psi_T\in L^2( \mathbb R^2_{x,y}
\times ]-\frac T2, \frac T2[)\,.
\end{equation}

We note that this functional has clearly a minimizer in $\mathcal
S_{\Omega}^{per}$ , where
\begin{equation}\label{Operunit}
\mathcal S_{\Omega}^{per}=\{ \Psi \in \mathcal
D_{\Omega}^{per}\,,\;
 \int_{\mathbb R^2_{x,y}\times ]-\frac T2, \frac T2[} |\Psi(x,y,z)|^2\; dx dy dz  =1\;\}.
\end{equation}
We denote this infimum by
\begin{equation}
E_\Omega^{per} = \inf_{\Psi \in \mathcal S_{\Omega}^{per}}
Q_{\Omega}^{per}(\Psi)\,.
\end{equation}
In the spirit of the Floquet theory,
 one can also be interested in
the analysis
 of the Floquet Bose-Einstein family of functionals, defined for $k
 \in \mathbb R$ by ~:
\begin{equation}\label{defQBEFlo}
\begin{array}{l}
\Psi \mapsto Q_{\Omega,k}^{Floq}(\Psi)\\
\quad := \int_{\mathbb R^2_{x,y}\times ]-\frac T2, \frac T2[}
   \left(\frac{1}{2}|(\nabla_{x,y} -i \Omega r_\perp) \Psi|^2 -\frac
     12  \Omega^2r^2|\Psi|^2\right.\\\left.
\quad\quad \quad 
+\frac 12 |(\partial_z
   + ik)\Psi|^2  + (V({\bf r}) +W_\ep(z)) |\Psi|^2\right)dxdydz\\
\quad\quad\quad\quad  + g \int_{\mathbb R^2_{x,y}\times ]-\frac T2, \frac T2[}|\Psi|^4\,dxdydz \,,
\end{array}
\end{equation}
where $\Psi$  satisfies
\begin{equation}
\Psi (x,y, z+T)= \Psi (x,y,z)\,.
\end{equation}
We call $\mathcal D_{\Omega,k}^{Floq}$, the natural maximal form domain
 of $ Q_{\Omega,k}^{Floq}$, which is actually independent of $k$ and $\Omega$~:
$$
\mathcal D_{\Omega,k}^{Floq}= \mathcal D_{\Omega}^{per}\,.
$$

We note that we have here a family over $k$ of functionals. Each of
these functionals has a minimizer in  $\mathcal
\mathcal S_{k}^{Floq}$, where
\begin{equation}\label{norma1b}
\mathcal S_{\Omega,k}^{Floq} =: \{ \Psi\in \mathcal D_{\Omega,k}^{Floq}
\,,\; \int_{\mathbb R^2_{x,y}\times ]-\frac T2, \frac T2[}
|\Psi|^2\; dx dy dz  =1\;\}.
\end{equation}
It is natural to be interested in the quantity
\begin{equation}\label{infkpsi}
E_\Omega^{Floq} : = \inf_k E_{\Omega,k}^{Floq}\,,
\end{equation}
with
\begin{equation}
E_{\Omega,k}^{Floq}:=\inf_{\Psi\in \mathcal S_{\Omega,k}^{Floq}}
Q_{\Omega,k}^{Floq}(\Psi)\,.
\end{equation}

\begin{remark}~\\
In the physics literature, this corresponds to a ground state of a
condensate at rest in the frame of the optical lattice. The energies
$E_{\Omega,k}^{Floq}$ describe states of the system where all the atoms
move with respect to the optical potential, at fixed velocity,
giving rise to  a constant current equal to $k$. Experimentally,
this is achieved by moving the  optical lattice with respect to the
condensate. We refer to \cite{KMPS} and the references therein.
\end{remark}

We would like to compare the various functionals.

Because $g\geq 0$, all the functional are  bounded from below, admit
an infimum and a lower bound is given by the analysis of the linear
 problem corresponding to $g=0$.
But there is probably no  existence of a minimizer for the Bose
Einstein functional. So we will have to consider
 minimizing sequences or approximate minimizers
 using possibly other models having minimizers.
 We refer for this to Subsection \ref{sspc}. The new point is
now that it is unclear that
 the various ground state energies obtained
 by the different procedures coincide~!

Finally, it could be interesting to compare directly the infimum of
 the three
 functionals above and then to find good approximations of these
 infima.

At the moment we have  shown in the ``linear'' section that~:
\begin{equation}
E_ \Omega(g=0)=E_\Omega^{per} (g=0) = E_\Omega^{Floq} (g=0)=\inf \sigma(H_\Omega)=\lambda_{1,z}+\omp \,.
\end{equation}
We have also mentioned, observing the monotonicity of the
 functionals
 with respect to $g$, that,  for $g\geq 0$,
\begin{equation}
\begin{array}{ll}
&E_\Omega(g=0)\leq E_\Omega(g) \,,\\
& E_\Omega^{per}(g=0) \leq
E_\Omega^{per}(g)\,,\\
\mbox{ and }&
 E_\Omega^{Floq} (g=0)\leq   E_\Omega^{Floq} (g)\,.
\end{array}
\end{equation}
This is further developed in Proposition \ref{pbsurR}.

\subsection{Comparison between the Floquet Bose-Einstein functionals
 and the periodic Bose-Einstein functional}\label{ssC1}
The argument which follows is only correct when $\Omega =0$. The
idea
 is that using the Kato inequality, one has
\begin{equation}
Q_{BE}^{per}(|\psi|)\leq Q_{k}^{Floq}(\psi)\,,\;\forall k
\end{equation}

In the periodic case, the minimization over complex functions leads
to the same infimum
 as in the real case.\\
We then obtain
\begin{equation}
E^{ Floq}\leq  E^{per} \leq E_{k}^{Floq}\,,
\end{equation}
hence
\begin{equation}
E^{Floq}= E^{per}\,.
\end{equation}
This seems to suggest that, when $\Omega=0$, there is no interest to
introduce Floquet conditions, if one is only interested
 in the determination of the ground state energy.
 \subsection{Comparison between the periodic Bose-Einstein functional
 and the   Bose-Einstein functional on $\mathbb R^3$}\label{ssC2}
Here we again work for general $\Omega$'s satisfying
\eqref{conditionsurOmega}.
We can naturally try what was working in the linear case. If
$\Psi^{per}(x,y,z)= \psi_{per}(x,y)\varphi_1(z)$
 is the periodic minimizer of the linear problem, i.e. satisfying
\begin{equation*}
H^\Omega \Psi^{per} = E^{lin}\Psi^{per}\,,\; ||\Psi^{per}||^2_{L^2(\mathbb
 R^2_{x,y}\times ]-\frac T2, \frac T2[)}=1\;
\end{equation*}
and
\begin{equation}\; \Psi^{per}(x,y,z+T) =\Psi^{per}(x,y,z)\,,
\end{equation}
with
$
 E^{lin}= \omp + \lambda_{1,z}$,
we can use as trial function for the Bose
 Einstein functional
\begin{equation}
\widetilde \Psi_R:= c_R\chi_R(z) \Psi^{per}(x,y,z)\,,
\end{equation}
where $c_R>0$ is determined by the condition
\begin{equation}
c_R^2\int_{\mathbb R^3}  | \chi_R(z) \Psi^{per}(x,y,z)|^2 \; dx dy
dz =1\,,
\end{equation}
Here $R\geq 1$ is a free parameter which will tend to $+\infty$
 and  $\chi_R$ is  a function with support in $[-RT-T, RT +T]$ equal to
$1$ on $[-RT,+RT]$ whose first derivative is independent of $R\geq 1$.\\
It is immediate to see, from the normalization chosen for $\Psi^{\rm
 per}$,  that
\begin{equation}
c_R \sim \frac{1}{\sqrt{2R}}\,,\; \mbox{ as } R\rightarrow +
\infty\,.
\end{equation}

We note that these trial functions are radial in the $(x,y)$ variable
 and are independent of $\Omega$.

The main point is to observe that  we have \beq
\begin{array}{ll}
\int_{\mathbb R^3} |\widetilde \Psi_R(x,y,z)|^4\; dxdydz &\sim 2R
c_R^4 \int_{\mathbb R^2\times]-\frac T2, \frac T2[} |
\Psi^{per}(x,y,z)|^4 \;dxdy dz \\ & \sim \frac{1}{2R} \int_{\mathbb
R^2\times]-\frac T2, \frac T2[} | \Psi^{per}(x,y,z)|^4 \;dxdy dz\,.
\end{array}
\eeq Hence
\begin{equation}
\lim_{R\rightarrow +\infty} \int_{\mathbb R^3} |\widetilde
\Psi_R(x,y,z)|^4\; dxdydz
 = 0\,,
\end{equation}\\
We also obtain easily that, as $R\ar +\infty$, \beq \lim_{R\ar
+\infty} \int_{\mathbb R^3}|(\nabla -i\Omega\times {\bf r}) \widetilde \Psi_R(x,y,z)|^2
 = \int_{\mathbb R^2\times ]-\frac T2, \frac
 T2[}|(\nabla-i\Omega\times {\bf r} )  \Psi^{per}(x,y,z)|^2\, dxdydz
\eeq
and
\begin{multline} \lim_{R\ar +\infty} \int_{\mathbb R^3} (
  V(x,y,z)-\frac 12
\Omega^2 r^2)|
\widetilde \Psi_R(x,y,z)|^2\\
 = \int_{\mathbb R^2\times ]-\frac T2,\frac T2[}(V(x,y,z)-\frac 12 \Omega^2
 r^2)  |\Psi^{per}(x,y,z)|^2\, dxdydz\,.
\end{multline}
 So we obtain that \beq
 \lim_{R\rightarrow +\infty}  Q_{\Omega}(\widetilde \Psi_R)=
E_\Omega^{per}(g=0)\,. \eeq

One can also observe that everything is actually $\Omega$ independent.\\

Combining with what we have verified for the linear case, we obtain
that,  for $g\geq 0$,
\begin{equation}
E ^{per}_\Omega(g=0)=E (g=0)\leq E_\Omega (g) \leq   \lim_{R\rightarrow +\infty}
Q_{\Omega}(\widetilde \Psi_R)= E_\Omega ^{per}(g=0)\,.
\end{equation}

So we have proved the
\begin{proposition}\label{pbsurR}~\\
For all $g\geq 0$, all $0\leq \Omega <\omp$,
\begin{equation}
E(g)= E_\Omega(g)= E _\Omega^{per} (g=0)= E_\Omega (g=0) =E(g=0)\,.
\end{equation}
\end{proposition}

The conclusion is that if we look at the Bose-Einstein functional on
$\mathbb R^3$ the infimum of the functional restricted to
$L^2$-normalized states is independent of $g\geq 0$  and $\Omega$ and is
immediately obtained
 by the ground state energy of the Hamiltonian attached to the case
 $g=0$
 and $\Omega=0$.\\

The result is of course also valid for the $1$-dimensional problem
and
is independent of any asymptotic analysis.\\
\begin{remark}~\\
Another natural physical problem would be to analyze the quantity
$$
\liminf_{N_c\ar +\infty}  \frac{1}{N_c} \left(
\inf_{\int_{-\frac{NT}{2}}^{+\frac {NT}{2}} |\Psi|^2\, dx = N_c }
Q_{\Omega}^{per,N}(\Psi)\right)
$$
or
$$
\limsup_{N_c\ar +\infty}  \frac{1}{N_c} \left(
\inf_{\int_{-\frac{NT}{2}}^{+\frac {NT}{2}} |\Psi|^2\, dx = N_c }
Q_{\Omega}^{per,N} (\Psi)\right)\,
$$
where we compute the energy by integrating over
 $N$  periods and where
 $$ N_c/N=\nu $$ ($\nu$ fixed). Upper bounds for this model
 are the periodic models with $g$ replaced by $g \nu$. We met already
 this problem in  Subsection
 \ref{lastsection}  for discrete models.
\end{remark}

\subsection{Comparison between the  $(NT)$-periodic problem and the
 $T$-problem}\label{sext}
In this subsection, we pursue the  analysis of
 the links between the $(NT)$-periodic
problem ($N>1$)  and the $T$-periodic problem.
We recall from Subsubsection \ref{sssNper} that, for the
$(NT)$-periodic problem in  $1D$,
 the ground state energy is $\lambda_{1,z}$. Moreover, in the
 semi-classical
 regime we have a packet of $N$ eigenvalues which are exponentially
 close
 separated  from the $(N+1)$-th eigenvalue by a splitting
 $\delta_z^{N}$ which satisfies \eqref{defphiN}. In case A, a natural
 question is~:\\
\begin{question}~\\
Is the minimizer of $\Epsilon_A^N$ $T$-periodic as in the linear
case ?
\end{question}
When the answer is yes, we immediately obtain that
\beq m_A^{N}
(\ep, \widehat g )= m_A (\ep,\frac{\widehat g}{N}) \eeq and we can
directly use what we have done for proving Theorem \ref{theoA}
  by replacing
$g$ by $\frac g N$.\\
To our knowledge the answer to this question is unknown, so it is
natural
 to look at simpler models.

In case B a natural question could be
\begin{question}~\\
Is the corresponding minimizer in this reduced space $T$-periodic ?
\end{question}
When the answer to this question is yes, we have seen in
\eqref{majcasB}
 together with the discussion around \eqref{tfN} (with the additional
 assumption that the minimizer is periodic)
  that
$$ m_{B,\Omega}^N =
 m_{B,\Omega}(\frac{\widetilde g}{N})\,.
$$
We could then use what we have used for the proof of Theorem
\ref{theoBWI}. This should work in the case when $g$ and $\Omega$
are small.

Actually, we could have asked more directly the following question~\\
\begin{question}~\\
Under which condition on $g$ and $\Omega$ is  the minimizer of the
$(NT)$-periodic initial problem $T$-periodic ?\end{question} If it
is the case, we get  immediately \beq \label{TNT} E_\Omega^{per,N} (g)=
E_\Omega^{per}(\frac{g}{N} ) \eeq So we can directly relate the treatment
of the $(NT)$-periodic problem to  our preceding studies, without
 any use of $(NT)$-periodic Wannier functions.\\

The general answer is unknown. One suspects
 by bifurcation arguments that it is true for $g$ and $\Omega$  small enough, but
 the physicists seem
 to wait for the other situation.
 This should in particular be  the case for
 sufficiently large
 rotation $\Omega$ (see for example \cite{CorR-DY1})
 or in the case with $(NT)$-Floquet conditions (see the discussion in
 Subsection \ref{lastsection} on
 an approximating model).


\appendix\section{Floquet theory}
 We follow the presentation of \cite{DiSj}
(p.~160-161), who are actually dealing with a more complicate
situation. If we take as lattice $\Gamma =\mathbb Z T $, the dual
lattice $\Gamma^*$
 is $\Gamma^*=\frac{2\pi}{T}  \mathbb Z$. For $u\in \mathcal S(\mathbb R)$
 and $k\in \mathbb R$, we put
\begin{equation}\label{A.1}
\mathcal U u(y;k) = \sum_{\gamma \in \Gamma} \exp i \gamma\,
k \,\; u(y-\gamma)\,.
\end{equation}
We notice that $\mathcal  U u(y;k)$ only depends on $k$ modulo
 the dual lattice $\Gamma^*$ so $\mathbb U u(y;k)$ is well
 defined
 on $\mathbb R \times \left(\mathbb R/\Gamma^*\right)$.\\

For $k \in \mathbb R/\Gamma^*$, we put
\begin{equation}\label{A.2}
\mathcal D'_k =\{u\in \mathcal D'(\mathbb R)\;;\;
u(y+\gamma)=\exp i\gamma \, k\,\; u(y)\}\,,
\end{equation}
and
\begin{equation}\label{A.3}
\mathcal H_k = \{u \in L^{2}_{loc}(\mathbb R)\cap \mathcal
D'_k\;;\;
 \int_{ E} |u(y)|^2\, dy <+\infty\}\,,
\end{equation}
where $E$ is a fundamental domain of $\Gamma$ (for example
$E=[-\frac T2, \frac T2[$).\\
When $k=0$, $\mathcal H_0$ denotes simply the space of
 $T$-periodic functions in $L^2_{loc}(\mathbb R)$.\\
~From \eqref{A.1} we see that
\begin{equation}\label{A.4}
\mathcal U u (\,\cdot\, ;k)\in \mathcal H_k\,,
\end{equation}
and if we view $\mathcal H_k$ as a bundle over $\mathbb
R/\Gamma^*$, we can view $\mathcal U u$ as a section of this
bundle. We write $\mathcal U u\in C^\infty(\mathbb
R/\Gamma^*;\mathcal H_k)$.\\

Now, if $v\in L^2(\mathbb R/\Gamma^*)$, one can expand it
 in a Fourier series~:
\begin{equation}\label{A.5}
v(k) = \sum_{\gamma\in \Gamma} \widehat v (\gamma)\, \exp i \gamma
\, k\,,
\end{equation}
where
\begin{equation}\label{A.6}
\widehat v (\gamma) = \frac {T}{2\pi} \int_{0}^{\frac {2\pi}{T}} \exp
-i\gamma\, k \,\; v(k)\; dk\,.
\end{equation}
We have the Parseval formula
\begin{equation}\label{A.7}
\frac {T}{2\pi} \int_0^{\frac{2\pi}{T}} |v(k)|^2\; dk =\sum_{\gamma
  \in \Gamma} |\widehat v (\gamma)|^2\,.
\end{equation}
Now, for any $y\in \mathbb R$, we can view \eqref{A.1} as the Fourier
expansion of $\mathcal U(y;\cdot)$, so \eqref{A.7} gives
\begin{equation}\label{A.8}
\frac {T}{2\pi} \int_0^{\frac{2\pi}{T}} |\mathcal U u (y;k)|^2\;dk
 = \sum_{\gamma \in \Gamma} |u(y-\gamma)|^2\,.
\end{equation}
Integrating over $y\in [-\frac T2, \frac T2]$, we get
\begin{equation}\label{A.9}
\frac {T}{2\pi} \int_{-\frac T2}^{\frac T2} \left(\int_0^{\frac{2\pi}{T}} |\mathcal U u
 (y;k)|^2\;dk\right) \;dy
 = \sum_{\gamma \in \Gamma}\int_{-\frac T2}^{\frac T2}  |u(y-\gamma)|^2\; dy =
 \int_{-\infty}^{+\infty} |u(y)|^2 \; dy\,.
\end{equation}
So $\mathcal U$ can be extended to an isometry from $L^2(\mathbb R)$
 into $L^2(\mathbb R/\Gamma^*;\mathcal H_k)$, which, as an
 Hilbert space,
 can be also described as the space of the $v$'s in $L^2_{loc}(\mathbb
 R^2)$
 such that
\begin{equation}\label{A.10}
\left\{
\begin{array}{ll}
v(y+ \gamma;k) = \exp i\gamma\, k\,\;
v(y;k)\,,&\forall \gamma \in \Gamma\,,\\
v(y;k + \gamma^*)= v(y;k)\,,&\forall \gamma^*\in \Gamma^*\,,
\end{array}\right.
\end{equation}
with the norm, whose square appears in the left hand side of
\eqref{A.9}, i.e.
\begin{equation}\label{A.11}
v \mapsto \sqrt{ \int_{[-\frac T2, \frac T2]\times [0,2\pi/T] }
  |v(y;k)|^2\;dy\,dk}\,.
\end{equation}
Like in the standard analysis of the Fourier transform, we have now to
analyze
 the surjectivity property and the construction of an inverse.\\
If $(y,k)\mapsto v(y;k)$
belongs to $C^\infty(\mathbb R/\Gamma^*;\mathcal H_k)$, we can
write the Fourier expansion of $v$ with respect to the second
variable~:
\begin{equation}\label{A.12}
v(y;k) = \sum_{\gamma\in \Gamma} \widehat v_\gamma (y)\,\exp
i\gamma\, k \,,
\end{equation}
with
\begin{equation}\label{A.13}
\widehat v_\gamma (y) = \frac{T}{2\pi} \int_{0}^{\frac{2\pi}{T}} \exp -i
\gamma\, k\,\; v(y;k) \;dk\,.
\end{equation}
Using the first line of \eqref{A.10}, we see that
\begin{equation}\label{A.14}
\widehat v_\gamma (y) =\widehat v_0 (y-\gamma)\,.
\end{equation}
Hence, when $v$ is smooth, we have
\begin{equation}\label{A.15}
v =\mathcal U \mathcal W v\,,
\end{equation}
where
\begin{equation}\label{A.16}
\mathcal W v(y)= \frac{T}{2\pi} \int_0^{\frac{2\pi}{T}} v(y;k)\;dk =
\widehat  v_0 (y)\,.
\end{equation}
Using \eqref{A.16}, \eqref{A.14} and \eqref{A.12}, we obtain
\begin{equation}\label{A.17}
\begin{array}{ll}
||\mathcal W v||^2_{L^2(\mathbb R)}
 &= \sum_{\gamma \in \Gamma} \int_{-\frac T2}^{\frac T2} |\widehat v_0(y-\gamma)|^2\;dy\\
&=  \sum_{\gamma \in \Gamma} \int_{-\frac T2}^{\frac T2} |\widehat v_\gamma(y)|^2\;dy\\
& = ||v||^2_{L^2(\mathbb R/\Gamma^*,\mathcal H_k)}\,.
\end{array}
\end{equation}
This shows that $\mathcal W$ is an isometry of  $L^2(\mathbb R/\Gamma^*;\mathcal H_k)$
 into
 $L^2(\mathbb R)$. Hence $\mathcal W$ is a bounded right inverse for
 $\mathcal U$. Since $\mathcal U$ is an isometry, we conclude that
 $\mathcal U$ is unitary and that its inverse is $\mathcal W$.\\
Now, when considering $P:=-\frac{d^2}{dy^2} + W(y)$ with $W$
 $T$-periodic, we see that
\begin{equation}\label{A.18}
\mathcal U P \mathcal U^{-1} = \int^{\oplus}_{[0,2\pi/T]}
\widetilde P_k \;dk\,,
\end{equation}
where,  by definition the right hand side in \eqref{A.18} denotes the
selfadjoint operator $Q$
 on  $L^2(\mathbb R/\Gamma^*;\mathcal H_k)$, with domain
 $L^2(\mathbb R/\Gamma^*;\mathcal H^2_k)$, which is  given by
\begin{equation}
Q v(y;k) = \left(
\widetilde P_k v(\cdot;k)\right)(y)\,.
\end{equation}
Here above
$$\mathcal H^2_k=\{u\in \mathcal H_k\,,\; u^{(\ell)} \in
\mathcal H_k\;\mbox{ for } |\ell|\leq 2\}\,,
$$
and $
\widetilde P_k$ is the selfadjoint operator on $\mathcal H_k$
associated to the
differential operator $-\frac{d^2}{dy^2} + W(y)$ with domain $\mathcal
H^2_k$.
\\
Finally, we note that for $k \in \mathbb R/\Gamma^*\,$, $
\widetilde P_k$
 is unitary equivalent to the operator
\begin{equation}\label{A.20}
 P_k = \exp -iy\, k\,\; \circ\,
\widetilde P_k \;\circ \,  \exp
iy\, k\,,
\end{equation}
which is now acting on $\mathcal H_0$ (corresponding to the
$T$-periodic fuunctions in $L^2_{loc}$).
More explicitly it takes the form
$$
 P_k = -\left(\frac{d}{dy} +i k\right)^2 + W(y)\,.
$$


\hfill

\noindent

{\bf Acknowledgements: } We are very grateful to Michiel Snoek
 for answering our questions on his PhDThesis and the physics knowledge 
  on the topic and to Xavier Blanc for his comments on the manuscript.
  We acknowledge support from the French ministry grant ANR-BLAN-0238, VoLQuan.


\bibliographystyle{plain}

\end{document}